\catcode`\@=11
\font\tensmc=cmcsc10      
\def\smc{\tensmc}

\def\hcorrection#1{\advance\hoffset by #1 }
\def\vcorrection#1{\advance\voffset by #1 }
\def\wlog#1{}
\newif\iftitle@
\outer\def\title{\title@true\vglue 24\p@ plus 12\p@ minus 12\p@
   \bgroup\let\\=\cr\tabskip\centering
   \halign to \hsize\bgroup\tenbf\hfill\ignorespaces##\unskip\hfill\cr}
\def\endtitle{\cr\egroup\egroup\vglue 18\p@ plus 12\p@ minus 6\p@}
\outer\def\author{\iftitle@\vglue -18\p@ plus -12\p@ minus -6\p@\fi\vglue
    12\p@ plus 6\p@ minus 3\p@\bgroup\let\\=\cr\tabskip\centering
    \halign to \hsize\bgroup\smc\hfill\ignorespaces##\unskip\hfill\cr}
\def\endauthor{\cr\egroup\egroup\vglue 18\p@ plus 12\p@ minus 6\p@}
\outer\def\heading{\bigbreak\bgroup\let\\=\cr\tabskip\centering
    \halign to \hsize\bgroup\smc\hfill\ignorespaces##\unskip\hfill\cr}
\def\endheading{\cr\egroup\egroup\nobreak\medskip}

\outer\def\proclaim#1{\medbreak\noindent\smc\ignorespaces
    #1\unskip.\enspace\sl\ignorespaces}
\outer\def\endproclaim{\par\ifdim\lastskip<\medskipamount\removelastskip
  \penalty 55 \fi\medskip\rm}
\outer\def\demo#1{\par\ifdim\lastskip<\smallskipamount\removelastskip
    \smallskip\fi\noindent{\smc\ignorespaces#1\unskip:\enspace}\rm
      \ignorespaces}
\outer\def\enddemo{\par\smallskip}
\newcount\footmarkcount@
\footmarkcount@=1
\def\makefootnote@#1#2{\insert\footins{\interlinepenalty=100
  \splittopskip=\ht\strutbox \splitmaxdepth=\dp\strutbox
  \floatingpenalty=\@MM
  \leftskip=\z@\rightskip=\z@\spaceskip=\z@\xspaceskip=\z@
  \noindent{#1}\footstrut\rm\ignorespaces #2\strut}}
\def\footnote{\let\@sf=\empty\ifhmode\edef\@sf{\spacefactor
   =\the\spacefactor}\/\fi\futurelet\next\footnote@}
\def\footnote@{\ifx"\next\let\next\footnote@@\else
    \let\next\footnote@@@\fi\next}
\def\footnote@@"#1"#2{#1\@sf\relax\makefootnote@{#1}{#2}}
\def\footnote@@@#1{$^{\number\footmarkcount@}$\makefootnote@
   {$^{\number\footmarkcount@}$}{#1}\global\advance\footmarkcount@ by 1 }

\hyphenation{man-u-script man-u-scripts ap-pen-dix ap-pen-di-ces}
\hyphenation{data-base data-bases}
\ifx\amstexloaded@\relax\catcode`\@=13
  \endinput\else\let\amstexloaded@=\relax\fi
\newlinechar=`\^^J
\def\eat@#1{}
\def\Space@.{\futurelet\Space@\relax}
\Space@. %
\newhelp\athelp@
{Only certain combinations beginning with @ make sense to me.^^J
Perhaps you wanted \string\@\space for a printed @?^^J
I've ignored the character or group after @.}
\def\futureletnextat@{\futurelet\next\at@}
{\catcode`\@=\active
\lccode`\Z=`\@ \lowercase
{\gdef@{\expandafter\csname futureletnextatZ\endcsname}
\expandafter\gdef\csname atZ\endcsname
   {\ifcat\noexpand\next a\def\next{\csname atZZ\endcsname}\else
   \ifcat\noexpand\next0\def\next{\csname atZZ\endcsname}\else
    \def\next{\csname atZZZ\endcsname}\fi\fi\next}
\expandafter\gdef\csname atZZ\endcsname#1{\expandafter
   \ifx\csname #1Zat\endcsname\relax\def\next
     {\errhelp\expandafter=\csname athelpZ\endcsname
      \errmessage{Invalid use of \string@}}\else
       \def\next{\csname #1Zat\endcsname}\fi\next}
\expandafter\gdef\csname atZZZ\endcsname#1{\errhelp
    \expandafter=\csname athelpZ\endcsname
      \errmessage{Invalid use of \string@}}}}
\def\atdef@#1{\expandafter\def\csname #1@at\endcsname}
\newhelp\defahelp@{If you typed \string\define\space cs instead of
\string\define\string\cs\space^^J
I've substituted an inaccessible control sequence so that your^^J
definition will be completed without mixing me up too badly.^^J
If you typed \string\define{\string\cs} the inaccessible control sequence^^J
was defined to be \string\cs, and the rest of your^^J
definition appears as input.}
\newhelp\defbhelp@{I've ignored your definition, because it might^^J
conflict with other uses that are important to me.}
\def\define{\futurelet\next\define@}
\def\define@{\ifcat\noexpand\next\relax
  \def\next{\define@@}%
  \else\errhelp=\defahelp@
  \errmessage{\string\define\space must be followed by a control
     sequence}\def\next{\def\garbage@}\fi\next}
\def\undefined@{}
\def\preloaded@{}
\def\define@@#1{\ifx#1\relax\errhelp=\defbhelp@
   \errmessage{\string#1\space is already defined}\def\next{\def\garbage@}%
   \else\expandafter\ifx\csname\expandafter\eat@\string
         #1@\endcsname\undefined@\errhelp=\defbhelp@
   \errmessage{\string#1\space can't be defined}\def\next{\def\garbage@}%
   \else\expandafter\ifx\csname\expandafter\eat@\string#1\endcsname\relax
     \def\next{\def#1}\else\errhelp=\defbhelp@
     \errmessage{\string#1\space is already defined}\def\next{\def\garbage@}%
      \fi\fi\fi\next}
\def\famzero{\fam\z@}

\def\lim{\mathop{\famzero lim}}

\def\textfont@#1#2{\def#1{\relax\ifmmode
    \errmessage{Use \string#1\space only in text}\else#2\fi}}
\textfont@\rm\tenrm
\textfont@\it\tenit
\textfont@\sl\tensl
\textfont@\bf\tenbf
\textfont@\smc\tensmc
\let\ic@=\/
\def\/{\unskip\ic@}
\def\textfonti{\the\textfont1 }
\def\t#1#2{{\edef\next{\the\font}\textfonti\accent"7F \next#1#2}}
\let\B=\=
\let\D=\.
\def~{\unskip\nobreak\ \ignorespaces}
{\catcode`\@=\active
\gdef\@{\char'100 }}
\atdef@-{\leavevmode\futurelet\next\athyph@}
\def\athyph@{\ifx\next-\let\next=\athyph@@
  \else\let\next=\athyph@@@\fi\next}
\def\athyph@@@{\hbox{-}}
\def\athyph@@#1{\futurelet\next\athyph@@@@}
\def\athyph@@@@{\if\next-\def\next##1{\hbox{---}}\else
    \def\next{\hbox{--}}\fi\next}
\def\.{.\spacefactor=\@m}
\atdef@.{\null.}
\atdef@,{\null,}
\atdef@;{\null;}
\atdef@:{\null:}
\atdef@?{\null?}
\atdef@!{\null!}
\def\srdr@{\thinspace}
\def\drsr@{\kern.02778em}
\def\sldl@{\kern.02778em}
\def\dlsl@{\thinspace}
\atdef@"{\unskip\futurelet\next\atqq@}
\def\atqq@{\ifx\next\Space@\def\next. {\atqq@@}\else
         \def\next.{\atqq@@}\fi\next.}
\def\atqq@@{\futurelet\next\atqq@@@}
\def\atqq@@@{\ifx\next`\def\next`{\atqql@}\else\def\next'{\atqqr@}\fi\next}
\def\atqql@{\futurelet\next\atqql@@}
\def\atqql@@{\ifx\next`\def\next`{\sldl@``}\else\def\next{\dlsl@`}\fi\next}
\def\atqqr@{\futurelet\next\atqqr@@}
\def\atqqr@@{\ifx\next'\def\next'{\srdr@''}\else\def\next{\drsr@'}\fi\next}

\def\textfontii{\the\textfont2 }
\def\{{\relax\ifmmode\lbrace\else
    {\textfontii f}\spacefactor=\@m\fi}
\def\}{\relax\ifmmode\rbrace\else
    \let\@sf=\empty\ifhmode\edef\@sf{\spacefactor=\the\spacefactor}\fi
      {\textfontii g}\@sf\relax\fi}
\def\nonhmodeerr@#1{\errmessage
     {\string#1\space allowed only within text}}
\def\linebreak{\relax\ifhmode\unskip\break\else
    \nonhmodeerr@\linebreak\fi}
\def\allowlinebreak{\relax
   \ifhmode\allowbreak\else\nonhmodeerr@\allowlinebreak\fi}
\newskip\saveskip@
\def\nolinebreak{\relax\ifhmode\saveskip@=\lastskip\unskip
  \nobreak\ifdim\saveskip@>\z@\hskip\saveskip@\fi
   \else\nonhmodeerr@\nolinebreak\fi}
\def\newline{\relax\ifhmode\null\hfil\break
    \else\nonhmodeerr@\newline\fi}
\def\nonmathaerr@#1{\errmessage
     {\string#1\space is not allowed in display math mode}}
\def\nonmathberr@#1{\errmessage{\string#1\space is allowed only in math mode}}
\def\mathbreak{\relax\ifmmode\ifinner\break\else
   \nonmathaerr@\mathbreak\fi\else\nonmathberr@\mathbreak\fi}
\def\nomathbreak{\relax\ifmmode\ifinner\nobreak\else
    \nonmathaerr@\nomathbreak\fi\else\nonmathberr@\nomathbreak\fi}
\def\allowmathbreak{\relax\ifmmode\ifinner\allowbreak\else
     \nonmathaerr@\allowmathbreak\fi\else\nonmathberr@\allowmathbreak\fi}
\def\pagebreak{\relax\ifmmode
   \ifinner\errmessage{\string\pagebreak\space
     not allowed in non-display math mode}\else\postdisplaypenalty-\@M\fi
   \else\ifvmode\penalty-\@M\else\edef\spacefactor@
       {\spacefactor=\the\spacefactor}\vadjust{\penalty-\@M}\spacefactor@
        \relax\fi\fi}
\def\nopagebreak{\relax\ifmmode
     \ifinner\errmessage{\string\nopagebreak\space
    not allowed in non-display math mode}\else\postdisplaypenalty\@M\fi
    \else\ifvmode\nobreak\else\edef\spacefactor@
        {\spacefactor=\the\spacefactor}\vadjust{\penalty\@M}\spacefactor@
         \relax\fi\fi}
\def\newpage{\relax\ifvmode\vfill\penalty-\@M\else\nonvmodeerr@\newpage\fi}
\def\nonvmodeerr@#1{\errmessage
    {\string#1\space is allowed only between paragraphs}}
\def\smallpagebreak{\relax\ifvmode\smallbreak
      \else\nonvmodeerr@\smallpagebreak\fi}
\def\medpagebreak{\relax\ifvmode\medbreak
       \else\nonvmodeerr@\medpagebreak\fi}
\def\bigpagebreak{\relax\ifvmode\bigbreak
      \else\nonvmodeerr@\bigpagebreak\fi}
\newdimen\captionwidth@
\captionwidth@=\hsize
\advance\captionwidth@ by -1.5in
\def\caption#1{}
\def\topspace#1{\gdef\thespace@{#1}\ifvmode\def\next
    {\futurelet\next\topspace@}\else\def\next{\nonvmodeerr@\topspace}\fi\next}
\def\topspace@{\ifx\next\Space@\def\next. {\futurelet\next\topspace@@}\else
     \def\next.{\futurelet\next\topspace@@}\fi\next.}
\def\topspace@@{\ifx\next\caption\let\next\topspace@@@\else
    \let\next\topspace@@@@\fi\next}
 \def\topspace@@@@{\topinsert\vbox to
       \thespace@{}\endinsert}
\def\topspace@@@\caption#1{\topinsert\vbox to
    \thespace@{}\nobreak
      \smallskip
    \setbox\z@=\hbox{\noindent\ignorespaces#1\unskip}%
   \ifdim\wd\z@>\captionwidth@
   \centerline{\vbox{\hsize=\captionwidth@\noindent\ignorespaces#1\unskip}}%
   \else\centerline{\box\z@}\fi\endinsert}
\def\midspace#1{\gdef\thespace@{#1}\ifvmode\def\next
    {\futurelet\next\midspace@}\else\def\next{\nonvmodeerr@\midspace}\fi\next}
\def\midspace@{\ifx\next\Space@\def\next. {\futurelet\next\midspace@@}\else
     \def\next.{\futurelet\next\midspace@@}\fi\next.}
\def\midspace@@{\ifx\next\caption\let\next\midspace@@@\else
    \let\next\midspace@@@@\fi\next}
 \def\midspace@@@@{\midinsert\vbox to
       \thespace@{}\endinsert}
\def\midspace@@@\caption#1{\midinsert\vbox to
    \thespace@{}\nobreak
      \smallskip
      \setbox\z@=\hbox{\noindent\ignorespaces#1\unskip}%
      \ifdim\wd\z@>\captionwidth@
    \centerline{\vbox{\hsize=\captionwidth@\noindent\ignorespaces#1\unskip}}%
    \else\centerline{\box\z@}\fi\endinsert}
\mathchardef\prime@="0230
\def\prime{{{}\prime@{}}}
\def\prim@s{\prime@\futurelet\next\pr@m@s}

\def\,{\relax\ifmmode\mskip\thinmuskip\else\thinspace\fi}
\def\!{\relax\ifmmode\mskip-\thinmuskip\else\negthinspace\fi}
\def\frac#1#2{{#1\over#2}}

\def\:{\nobreak\hskip.1111em{:}\hskip.3333em plus .0555em\relax}
\def\intic@{\mathchoice{\hskip5\p@}{\hskip4\p@}{\hskip4\p@}{\hskip4\p@}}
\def\negintic@
 {\mathchoice{\hskip-5\p@}{\hskip-4\p@}{\hskip-4\p@}{\hskip-4\p@}}
\def\intkern@{\mathchoice{\!\!\!}{\!\!}{\!\!}{\!\!}}
\def\intdots@{\mathchoice{\cdots}{{\cdotp}\mkern1.5mu
    {\cdotp}\mkern1.5mu{\cdotp}}{{\cdotp}\mkern1mu{\cdotp}\mkern1mu
      {\cdotp}}{{\cdotp}\mkern1mu{\cdotp}\mkern1mu{\cdotp}}}
\newcount\intno@
\def\iint{\intno@=\tw@\futurelet\next\ints@}
\def\iiint{\intno@=\thr@@\futurelet\next\ints@}
\def\iiiint{\intno@=4 \futurelet\next\ints@}
\def\idotsint{\intno@=\z@\futurelet\next\ints@}
\def\ints@{\findlimits@\ints@@}
\newif\iflimtoken@
\newif\iflimits@
\def\findlimits@{\limtoken@false\limits@false\ifx\next\limits
 \limtoken@true\limits@true\else\ifx\next\nolimits\limtoken@true\limits@false
    \fi\fi}
\def\multintlimits@{\intop\ifnum\intno@=\z@\intdots@
  \else\intkern@\fi
    \ifnum\intno@>\tw@\intop\intkern@\fi
     \ifnum\intno@>\thr@@\intop\intkern@\fi\intop}
\def\multint@{\int\ifnum\intno@=\z@\intdots@\else\intkern@\fi
   \ifnum\intno@>\tw@\int\intkern@\fi
    \ifnum\intno@>\thr@@\int\intkern@\fi\int}
\def\ints@@{\iflimtoken@\def\ints@@@{\iflimits@
   \negintic@\mathop{\intic@\multintlimits@}\limits\else
    \multint@\nolimits\fi\eat@}\else
     \def\ints@@@{\multint@\nolimits}\fi\ints@@@}
\def\Sb{_\bgroup\vspace@
        \baselineskip=\fontdimen10 \scriptfont\tw@
        \advance\baselineskip by \fontdimen12 \scriptfont\tw@
        \lineskip=\thr@@\fontdimen8 \scriptfont\thr@@
        \lineskiplimit=\thr@@\fontdimen8 \scriptfont\thr@@
        \Let@\vbox\bgroup\halign\bgroup \hfil$\scriptstyle
            {##}$\hfil\cr}
\def\endSb{\crcr\egroup\egroup\egroup}
\def\Sp{^\bgroup\vspace@
        \baselineskip=\fontdimen10 \scriptfont\tw@
        \advance\baselineskip by \fontdimen12 \scriptfont\tw@
        \lineskip=\thr@@\fontdimen8 \scriptfont\thr@@
        \lineskiplimit=\thr@@\fontdimen8 \scriptfont\thr@@
        \Let@\vbox\bgroup\halign\bgroup \hfil$\scriptstyle
            {##}$\hfil\cr}
\def\endSp{\crcr\egroup\egroup\egroup}
\def\Let@{\relax\iffalse{\fi\let\\=\cr\iffalse}\fi}
\def\vspace@{\def\vspace##1{\noalign{\vskip##1 }}}
\def\aligned{\,\vcenter\bgroup\vspace@\Let@\openup\jot\m@th\ialign
  \bgroup \strut\hfil$\displaystyle{##}$&$\displaystyle{{}##}$\hfil\crcr}
\def\endaligned{\crcr\egroup\egroup}
\def\matrix{\,\vcenter\bgroup\Let@\vspace@
    \normalbaselines
  \m@th\ialign\bgroup\hfil$##$\hfil&&\quad\hfil$##$\hfil\crcr
    \mathstrut\crcr\noalign{\kern-\baselineskip}}
\def\endmatrix{\crcr\mathstrut\crcr\noalign{\kern-\baselineskip}\egroup
                \egroup\,}
\newtoks\hashtoks@
\hashtoks@={#}
\def\format{\crcr\egroup\iffalse{\fi\ifnum`}=0 \fi\format@}
\def\format@#1\\{\def\preamble@{#1}%
  \def\c{\hfil$\the\hashtoks@$\hfil}%
  \def\r{\hfil$\the\hashtoks@$}%
  \def\l{$\the\hashtoks@$\hfil}%
  \setbox\z@=\hbox{\xdef\Preamble@{\preamble@}}\ifnum`{=0 \fi\iffalse}\fi
   \ialign\bgroup\span\Preamble@\crcr}
 
\let\hdots=\ldots
\def\cases{\left\{\,\vcenter\bgroup\vspace@
     \normalbaselines\openup\jot\m@th
       \Let@\ialign\bgroup$##$\hfil&\quad$##$\hfil\crcr
      \mathstrut\crcr\noalign{\kern-\baselineskip}}

\newif\iftagsleft@
\tagsleft@true
\def\TagsOnRight{\global\tagsleft@false}
\def\tag#1$${\iftagsleft@\leqno\else\eqno\fi
 \hbox{\def\pagebreak{\global\postdisplaypenalty-\@M}%
 \def\nopagebreak{\global\postdisplaypenalty\@M}\rm(#1\unskip)}%
  $$\postdisplaypenalty\z@\ignorespaces}
\interdisplaylinepenalty=\@M
\def\allowdisplaybreak@{\def\allowdisplaybreak{\noalign{\allowbreak}}}
\def\displaybreak@{\def\displaybreak{\noalign{\break}}}
\def\align#1\endalign{\def\tag{&}\vspace@\allowdisplaybreak@\displaybreak@
  \iftagsleft@\lalign@#1\endalign\else
   \ralign@#1\endalign\fi}
\def\ralign@#1\endalign{\displ@y\Let@\tabskip\centering\halign to\displaywidth
     {\hfil$\displaystyle{##}$\tabskip=\z@&$\displaystyle{{}##}$\hfil
       \tabskip=\centering&\llap{\hbox{(\rm##\unskip)}}\tabskip\z@\crcr
             #1\crcr}}
\def\lalign@
 #1\endalign{\displ@y\Let@\tabskip\centering\halign to \displaywidth
   {\hfil$\displaystyle{##}$\tabskip=\z@&$\displaystyle{{}##}$\hfil
   \tabskip=\centering&\kern-\displaywidth
        \rlap{\hbox{(\rm##\unskip)}}\tabskip=\displaywidth\crcr
               #1\crcr}}
\def\overrightarrow{\mathpalette\overrightarrow@}
\def\overrightarrow@#1#2{\vbox{\ialign{$##$\cr
    #1{-}\mkern-6mu\cleaders\hbox{$#1\mkern-2mu{-}\mkern-2mu$}\hfill
     \mkern-6mu{\to}\cr
     \noalign{\kern -1\p@\nointerlineskip}
     \hfil#1#2\hfil\cr}}}
\def\overleftarrow{\mathpalette\overleftarrow@}
\def\overleftarrow@#1#2{\vbox{\ialign{$##$\cr
     #1{\leftarrow}\mkern-6mu\cleaders\hbox{$#1\mkern-2mu{-}\mkern-2mu$}\hfill
      \mkern-6mu{-}\cr
     \noalign{\kern -1\p@\nointerlineskip}
     \hfil#1#2\hfil\cr}}}
\def\overleftrightarrow{\mathpalette\overleftrightarrow@}
\def\overleftrightarrow@#1#2{\vbox{\ialign{$##$\cr
     #1{\leftarrow}\mkern-6mu\cleaders\hbox{$#1\mkern-2mu{-}\mkern-2mu$}\hfill
       \mkern-6mu{\to}\cr
    \noalign{\kern -1\p@\nointerlineskip}
      \hfil#1#2\hfil\cr}}}
\def\underrightarrow{\mathpalette\underrightarrow@}
\def\underrightarrow@#1#2{\vtop{\ialign{$##$\cr
    \hfil#1#2\hfil\cr
     \noalign{\kern -1\p@\nointerlineskip}
    #1{-}\mkern-6mu\cleaders\hbox{$#1\mkern-2mu{-}\mkern-2mu$}\hfill
     \mkern-6mu{\to}\cr}}}
\def\underleftarrow{\mathpalette\underleftarrow@}
\def\underleftarrow@#1#2{\vtop{\ialign{$##$\cr
     \hfil#1#2\hfil\cr
     \noalign{\kern -1\p@\nointerlineskip}
     #1{\leftarrow}\mkern-6mu\cleaders\hbox{$#1\mkern-2mu{-}\mkern-2mu$}\hfill
      \mkern-6mu{-}\cr}}}
\def\underleftrightarrow{\mathpalette\underleftrightarrow@}
\def\underleftrightarrow@#1#2{\vtop{\ialign{$##$\cr
      \hfil#1#2\hfil\cr
    \noalign{\kern -1\p@\nointerlineskip}
     #1{\leftarrow}\mkern-6mu\cleaders\hbox{$#1\mkern-2mu{-}\mkern-2mu$}\hfill
       \mkern-6mu{\to}\cr}}}
\def\sqrt#1{\radical"270370 {#1}}
\def\dots{\relax\ifmmode\let\next=\ldots\else\let\next=\tdots@\fi\next}
\def\tdots@{\unskip\ \tdots@@}
\def\tdots@@{\futurelet\next\tdots@@@}
\def\tdots@@@{$\mathinner{\ldotp\ldotp\ldotp}\,
   \ifx\next,$\else
   \ifx\next.\,$\else
   \ifx\next;\,$\else
   \ifx\next:\,$\else
   \ifx\next?\,$\else
   \ifx\next!\,$\else
   $ \fi\fi\fi\fi\fi\fi}
\def\text{\relax\ifmmode\let\next=\text@\else\let\next=\text@@\fi\next}
\def\text@@#1{\hbox{#1}}
\def\text@#1{\mathchoice
 {\hbox{\everymath{\displaystyle}\def\textfonti{\the\textfont1 }%
    \def\textfontii{\the\textfont2 }\textdef@@ T#1}}
 {\hbox{\everymath{\textstyle}\def\textfonti{\the\textfont1 }%
    \def\textfontii{\the\textfont2 }\textdef@@ T#1}}
 {\hbox{\everymath{\scriptstyle}\def\textfonti{\the\scriptfont1 }%
   \def\textfontii{\the\scriptfont2 }\textdef@@ S\rm#1}}
 {\hbox{\everymath{\scriptscriptstyle}\def\textfonti{\the\scriptscriptfont1 }%
   \def\textfontii{\the\scriptscriptfont2 }\textdef@@ s\rm#1}}}
\def\textdef@@#1{\textdef@#1\rm \textdef@#1\bf
   \textdef@#1\sl \textdef@#1\it}

\def\textdef@#1#2{\def\next{\csname\expandafter\eat@\string#2fam\endcsname}%
\if S#1\edef#2{\the\scriptfont\next\relax}%
 \else\if s#1\edef#2{\the\scriptscriptfont\next\relax}%
 \else\edef#2{\the\textfont\next\relax}\fi\fi}
\scriptfont\itfam=\tenit \scriptscriptfont\itfam=\tenit
\scriptfont\slfam=\tensl \scriptscriptfont\slfam=\tensl
\mathcode`\0="0030
\mathcode`\1="0031
\mathcode`\2="0032
\mathcode`\3="0033
\mathcode`\4="0034
\mathcode`\5="0035
\mathcode`\6="0036
\mathcode`\7="0037
\mathcode`\8="0038
\mathcode`\9="0039
\def\Cal{\relax\ifmmode\let\next=\Cal@\else
     \def\next{\errmessage{Use \string\Cal\space only in math mode}}\fi\next}
\def\Cal@#1{{\fam2 #1}}
\def\bold{\relax\ifmmode\let\next=\bold@\else
   \def\next{\errmessage{Use \string\bold\space only in math
      mode}}\fi\next}\def\bold@#1{{\fam\bffam #1}}
\mathchardef\Gamma="0000
\mathchardef\Delta="0001
\mathchardef\Theta="0002
\mathchardef\Lambda="0003
\mathchardef\Xi="0004
\mathchardef\Pi="0005
\mathchardef\Sigma="0006
\mathchardef\Upsilon="0007
\mathchardef\Phi="0008
\mathchardef\Psi="0009
\mathchardef\Omega="000A
\mathchardef\varGamma="0100
\mathchardef\varDelta="0101
\mathchardef\varTheta="0102
\mathchardef\varLambda="0103
\mathchardef\varXi="0104
\mathchardef\varPi="0105
\mathchardef\varSigma="0106
\mathchardef\varUpsilon="0107
\mathchardef\varPhi="0108
\mathchardef\varPsi="0109
\mathchardef\varOmega="010A
\font\dummyft@=dummy
\fontdimen1 \dummyft@=\z@
\fontdimen2 \dummyft@=\z@
\fontdimen3 \dummyft@=\z@
\fontdimen4 \dummyft@=\z@
\fontdimen5 \dummyft@=\z@
\fontdimen6 \dummyft@=\z@
\fontdimen7 \dummyft@=\z@
\fontdimen8 \dummyft@=\z@
\fontdimen9 \dummyft@=\z@
\fontdimen10 \dummyft@=\z@
\fontdimen11 \dummyft@=\z@
\fontdimen12 \dummyft@=\z@
\fontdimen13 \dummyft@=\z@
\fontdimen14 \dummyft@=\z@
\fontdimen15 \dummyft@=\z@
\fontdimen16 \dummyft@=\z@
\fontdimen17 \dummyft@=\z@
\fontdimen18 \dummyft@=\z@
\fontdimen19 \dummyft@=\z@
\fontdimen20 \dummyft@=\z@
\fontdimen21 \dummyft@=\z@
\fontdimen22 \dummyft@=\z@
\def\fontlist@{\\{\tenrm}\\{\sevenrm}\\{\fiverm}\\{\teni}\\{\seveni}%
 \\{\fivei}\\{\tensy}\\{\sevensy}\\{\fivesy}\\{\tenex}\\{\tenbf}\\{\sevenbf}%
 \\{\fivebf}\\{\tensl}\\{\tenit}\\{\tensmc}}
\def\dodummy@{{\def\\##1{\global\let##1=\dummyft@}\fontlist@}}
\newif\ifsyntax@
\newcount\countxviii@
\def\newtoks@{\alloc@5\toks\toksdef\@cclvi}
\def\nopages@{\output={\setbox\z@=\box\@cclv \deadcycles=\z@}\newtoks@\output}
\def\syntax{\syntax@true\dodummy@\countxviii@=\count18
\loop \ifnum\countxviii@ > \z@ \textfont\countxviii@=\dummyft@
   \scriptfont\countxviii@=\dummyft@ \scriptscriptfont\countxviii@=\dummyft@
     \advance\countxviii@ by-\@ne\repeat
\dummyft@\tracinglostchars=\z@
  \nopages@\frenchspacing\hbadness=\@M}
\def\magstep#1{\ifcase#1 1000\or
 1200\or 1440\or 1728\or 2074\or 2488\or
 \errmessage{\string\magstep\space only works up to 5}\fi\relax}
{\lccode`\2=`\p \lccode`\3=`\t
 \lowercase{\gdef\tru@#123{#1truept}}}

\def\scaletype#1{\mag=#1\relax
 \hsize=\expandafter\tru@\the\hsize
 \vsize=\expandafter\tru@\the\vsize
 \dimen\footins=\expandafter\tru@\the\dimen\footins}

\def\scalefont#1#2\andcallit#3{\edef\font@{\the\font}#1\font#3=
  \fontname\font\space scaled #2\relax\font@}
\def\Mag@#1#2{\ifdim#1<1pt\multiply#1 #2\relax\divide#1 1000 \else
  \ifdim#1<10pt\divide#1 10 \multiply#1 #2\relax\divide#1 100\else
  \divide#1 100 \multiply#1 #2\relax\divide#1 10 \fi\fi}
\def\scalelinespacing#1{\Mag@\baselineskip{#1}\Mag@\lineskip{#1}%
  \Mag@\lineskiplimit{#1}}
\def\wlog#1{\immediate\write-1{#1}}
\catcode`\@=\active

\input vanilla.sty
\magnification \magstep 1
\baselineskip = 12pt
\voffset = .2cm
\vglue 1.2cm
\hsize 13.34cm
\vsize 20.14cm
\parindent = 0mm
\font\smrm = cmr9
\font\big =cmb10 scaled \magstep 2
\TagsOnRight

\title \big
On Lunn-Senior's Mathematical Model \\
\big of Isomerism in Organic Chemistry. Part II

\endtitle

\author
Valentin Vankov Iliev
\footnote"*"{\smrm Partially supported by Grant MM-1106/2001 of the
Bulgarian Foundation of Scientific Rese\-arch}

\endauthor

\centerline {\it Section of Algebra, Institute of Mathematics and Informatics}
\centerline {\it Bulgarian Academy of Sciences, 1113 Sofia, Bulgaria}
\centerline {\it E-mail: viliev\@math.bas.bg, viliev\@aubg.bg}

\heading
7. Introduction

\endheading

This paper is a continuation of [3], and we shall use all terminology and
notation introduced there.

7.1.  Let a molecule's skeleton with $d$ unsatisfied single valences is fixed.
Lunn-Senior's thesis 1.5.1 asserts that for any type of isomerism of the
molecule with that skeleton (univalent substitution isomerism, stereoisomerism,
or structural isomerism) there exists a permutation group $W\leq S_d$ such that
each isomer of that type can be identified with a $W$-orbit in the set $T_d$,
that is, with an element of the set $T_{d;W}$ of $W$-orbits.  If $W=G$, where
$G$ is Lunn-Senior's group of substitution isomerism, then the members of the
set $T_{d;G}$ represent the corresponding univalent substitution isomers, and
the inequalities $a<b$, $a,b\in T_{d;G}$, represent the substitution reactions
$b\longrightarrow a$ among them.  We denote by $G^{\prime}$ and
$G^{\prime\prime}$ the groups of stereoisomerism and structural isomerism,
respectively, of the molecule under question.  The group $G^\prime$ either
coincides with $G$ (there are no chiral pairs), or $G\leq G^\prime$ and
$|G^\prime :G|=2$ (otherwise) --- see [7, V] or 1.5.1.

The aim of the paper is to consider a mathematical model in which the following
statement from chemistry becomes a well formed formula:  the derivatives that
correspond to the structural formulae $a,b\in T_{d;G}$ can not be distinguished
via substitution reactions.  The present approach gives a conceptual basis of
the Lunn-Senior's {\it ad hoc} considerations in [7, VI] (see also 6.4) on
the diamers of ethene and of our statements in [5, Sections 3,4,5] concerning
the identification of the isomers of certain compounds with one
mono-substitution and at least three di-substitution homogeneous derivatives
(for instance, benzene and cyclopropane).

It is natural to suppose that $a$ and $b$ are identical as structural isomers,
that is, $a$ and $b$ are contained in one and the same
$G^{\prime\prime}$-orbit, --- otherwise they are trivially distinguishable;
thus they also have the same empirical formula:  $a,b\in T_{\lambda;G}$,
$\lambda\in P_d$.  Since the substitution reactions are modeled by the
inequalities in $T_{d;G}$, and since the symmetries are the automorphisms of
the corresponding algebraic structure --- the partially ordered set $T_{d;G}$
--- we define that the products corresponding to $a$ and $b$ are {\it
indistinguishable via substitution reactions} if $a$ and $b$ are identical as
structural isomers and if there exists an automorphism $\alpha\colon T_{d;G}\to
T_{d;G}$ of the partially ordered set $T_{d;G}$, such that:  (a)
$\alpha(T_{\mu;G})=T_{\mu;G}$ for any $\mu\in P_d$, (b) $\alpha$ maps any
chiral pair onto a chiral pair, and (c) $\alpha(a)=b$.  Otherwise, they are
called {\it distinguishable via substitution reactions}.  The automorphisms
$\alpha$ of $T_{d;G}$, which satisfy condition (a) form a group
$Aut_0(T_{d;G})$ that acts naturally on the set $T_{d;G}$.  The elements
$\alpha\in Aut_0(T_{d;G})$ that obey (b) form a subgroup $Aut_0^\prime
(T_{d;G})\leq Aut_0(T_{d;G})$.

In practice, however, the chemists know experimentally enough derivatives as
well as enough substitution reactions among them only for a small number of
$\lambda\in P_d$ (mono-substitution, di-substitution, tri-substitution
homogeneous, etc., derivatives). Thus, we are forced to consider subsets of
$T_{d;G}$ consisting of several $T_{\mu;G}$, $\mu\in P_d$, with the
induced partial order.  For any subset $D\subset P_d$ we define the partially
ordered set $T_{D;G}=\cup_{\mu\in D} T_{\mu;G}$ and generalize the
above definition in the following natural way:  the products corresponding to
$a, b\in T_{\lambda;G}$, $\lambda\in D$, are {\it indistinguishable via
substitution reactions among the elements of $T_{D;G}$} if $a$ and $b$ are
identical as structural isomers, and if there exists an automorphism
$\alpha\colon T_{D;G}\to T_{D;G}$ of the partially ordered set $T_{D;G}$, such
that (a$_D$) $\alpha(T_{\mu;G})=T_{\mu;G}$ for any $\mu\in D$, (b$_D$)
$\alpha$ maps any chiral pair onto a chiral pair, and (c$_D$) $\alpha(a)=b$.
Otherwise, they are said to be {\it distinguishable via substitution reactions
among the elements of $T_{D;G}$}.  The automorphisms $\alpha$ of $T_{D;G}$
possessing property (a$_D$) form a group $Aut_0(T_{D;G})$, and those
$\alpha\in Aut_0(T_{D;G})$ which obey (b$_D$) form a subgroup
$Aut_0^\prime(T_{D;G})\leq Aut_0(T_{D;G})$.  The elements of the group
$Aut_0^\prime(T_{D;G})$ could be called {\it chiral automorphisms} of the
partially ordered set $T_{D;G}$.

Obviously, if $D\subset D^\prime$, then any automorphism from the group
$Aut_0^\prime (T_{D^\prime;G})$ induces by restriction an automorphism from the
group $Aut_0^\prime (T_{D;G})$, and we obtain a homomorphism of groups
$Aut_0^\prime (T_{D^\prime;G})\to Aut_0^\prime (T_{D;G})$.  The fact that this
homomorphism is, in general, not surjective reflects the plain observation that
the products corresponding to $a, b\in T_{\lambda;G}$, $\lambda\in D$, can be
indistinguishable via substitution reactions among the elements of $T_{D;G}$,
but distinguishable via substitution reactions among the elements of (the wider
set) $T_{D^\prime;G}$.

7.2.  In the beginning of Section 8 we study the relation between the partially
ordered sets $T_{d;W}$ and $T_{d;W^\prime}$, for a couple of groups $W\leq
W^\prime\leq S_d$ with $|W^\prime :W|=2$, and the formal behaviour of the group
$Aut_0^{W^\prime}(T_{D;W})$ of chiral automorphisms of $T_{d;W}$.  Thus, when
$W=G$, $W^\prime=G^\prime$ we have $Aut_0^{W^\prime}(T_{D;W})=
Aut_0^\prime(T_{D;W})$, and we generalize the \lq\lq chiral" situation.

The renumberings of the $d$ unsatisfied valences of the skeleton, which
do not affect the group $G$ of symmetries of the molecule, are exactly the
elements $\nu$ of the normalizer $N$ of $G$ in $S_d$, and any such renumbering
$\nu$ induces an automorphism $\hat{\nu}\in Aut_0(T_{D;G})$, $D\subset P_d$.
In addition, if $\nu\in N^\prime$, where $N^\prime$ is the intersection of $N$
with the normalizer of $G^\prime$ in $S_d$, then $\hat{\nu}\in
Aut_0^\prime(T_{D;G})$.  After A.  Kerber (see [6, 1.5]), we call the
automorphisms $\hat{\nu}$, where $\nu\in N^\prime$,  {\it hidden symmetries}
of the molecule under consideration, or {\it hidden automorphisms} of
the partially ordered set $T_{D;W}$.

In the formal treatment, Theorem 8.2.2 yields
$\hat{\nu}\in Aut_0^{W^\prime}(T_{D;W})$, $D\subset P_d$, and, in particular,
shows that the hidden symmetries form a subgroup $Aut_0^{h;W^\prime}(T_{D;W})$
of $Aut_0^{W^\prime}(T_{D;W})$, isomorphic to the factor-group $N^\prime/W$
when $(1^d)\in D$.  In case $W=G$, $W^\prime=G^\prime$ we denote
$Aut_0^{h;W^\prime}(T_{D;W})$ by $Aut_0^h(T_{D;G})$.  Under certain
condition, chemically equivalent to existence of chiral pairs among the
derivatives with empirical formulae from $D$, the group
$Aut_0^{h;W^\prime}(T_{D;W})$ contains a copy of the factor-group $W^\prime
/W$, generated by the so called {\it chiral involution} $\hat\tau$ that
permutes the members of any chiral pair and leaves invariant the diamers
(Corollary 8.2.4).  In particular, the members of a chiral pair are
indistinguishable via substitution reactions, which is the content of Theorem
8.2.6.

On the other hand, the normalizer $N$ acts in a natural way on the set $X_W$ of
all one-dimensional characters $\chi$ of the group $W$, $\chi\mapsto\nu\chi$,
and this action induces actions of both the factor-group $N/W$ and, under some
condition on $D$, of the group $Aut_0^{h;W^\prime}(T_{D;W})$ of hidden
automorphisms, on the set $X_W$.  Given $\theta\in X_{S_\lambda}$, $\lambda\in
P_d$, any automorphism $\hat{\nu}$, $\nu\in N$, maps the set of all
$(\chi,\theta)$-orbits onto the set of all $(\nu\chi,\theta)$-orbits in $T_d$,
and, in particular, any hidden symmetry $\alpha\in Aut_0^{h;W^\prime}(T_{D;W})$
does the same.  This is proved in Theorem 9.1.1.  As Corollary 9.1.2 we get
that any automorphism $\hat{\nu}$, $\nu\in N$, transforms the set of all
$\chi$-orbits onto the set of all $\nu\chi$-orbits.

The Extended Lunn-Senior's thesis 1.6.1, hypothesis 5, asserts that for any
$\lambda\in P_d$, and for any pair of characters $(\chi,\theta)\in X_W\times
X_{S_\lambda}$ the subset $T_{\lambda;\chi,\theta}$ of $T_{\lambda;W}$
represents a type property $(\chi,\theta)$ of the molecule.  Moreover, the
automorphism $\hat\nu$, $\nu\in N^\prime$, of our algebraic structure
transforms the set $T_{\lambda;\chi,\theta}$ onto the set
$T_{\lambda;\nu\chi.\theta}$, so the type properties $(\chi,\theta)$ and
$(\nu\chi,\theta)$ are not distinguishable.  Thus, it is natural to define that
the derivatives which correspond to the structural formulae $a,b\in
T_{\lambda;G}$, $\lambda\in P_d$, are {\it indistinguishable via pairs of
characters} if $a$ and $b$ are structurally identical, if for any pair of
characters $(\chi,\theta)\in X_G\times X_{S_\lambda}$ with $a\in
T_{\lambda;\chi,\theta}$ and $b\notin T_{\lambda;\chi,\theta}$, there exists a
$\nu\in N^\prime$ such that $b\in T_{\lambda;\nu\chi,\theta}$ and $a\notin
T_{\lambda;\nu\chi,\theta}$, and if for any pair of characters
$(\chi,\theta)\in X_G\times X_{S_\lambda}$ with $b\in T_{\lambda;\chi,\theta}$
and $a\notin T_{\lambda;\chi,\theta}$, there exists a $\nu\in N^\prime$ such
that $a\in T_{\lambda;\nu\chi,\theta}$ and $b\notin
T_{\lambda;\nu\chi,\theta}$.  By the specialization $\theta=1_{S_\lambda}$, we
obtain the definition of indistinguishability via characters.  Otherwise, these
derivatives are called {\it distinguishable via pairs of characters},
respectively, {\it via characters}.  The central results here are Theorem 9.2.1
and Corollary 9.2.2, which assert that the members of a chiral pair are
indistinguishable via pairs of characters and via characters.

In Section 10 we illustrate these ideas by several molecules:  those ethene
$C_2H_4$, benzene $C_6H_6$, and cyclopropane $C_3H_6$.  The results on the last
two are valid in a more general situation described in [5].  Since there are no
chiral pairs among the derivatives of ethene, its group of \lq\lq chiral"
automorphisms is relatively big --- it is described in Theorem 10.1.1.  In
Proposition 10.1.2 the group of the hidden symmetries of ethene is described.
Corollary 10.1.3 gives a rigorous treatment of Lunn-Senior's considerations in
[7, VI], concerning the indistiguishability of the diameric pairs
$\{a_{\left(2^2\right)},b_{\left(2^2\right)}\}$,
$\{a_{\left(2,1^2\right)},b_{\left(2,1^2\right)}\}$ of ethene.  Proposition
10.1.4 and Corollary 10.1.5 specify their indistiguishability via (pairs of )
characters, mentioned in 6.4.  The group of chiral automorphisms
$Aut_0^\prime(T_{D;G})$ of benzene for $D=\{(4,2),(3^2)\}$ is the unit group,
so we get once again that the K\"orner relations identify completely the
$(4,2)$- and $(3^2)$-derivatives of benzene --- see Theorem 10.2.1 and
Corollary 10.2.2.  Theorem 10.3.5 asserts that the group of chiral
automorphisms $Aut_0^\prime(T_{D;G})$ of cyclopropane for
$D=\{(6),(5,1),(4,2),(4,1^2),(3^2)\}$ is direct product of two cyclic groups of
order $2$.  Corollary 10.3.6 yields that for this $D$ all chiral pairs and all
diamers are pairwise distinguishable via substitution reactions.

\heading
8. Automorphisms

\endheading

8.1.  Let $W\leq S_d$ be a permutation group.  We set $T_D=\cup_{\mu\in
D}T_\mu$, $T_{D;W}=\cup_{\mu\in D}T_{\mu;W}$, and consider the partial order on
$T_{D;W}$, induced from $T_{d;W}$.  By definition, a map $u\colon T_{D;W}\to
T_{D;W}$ is an automorphism of the partially ordered set $T_{D;W}$ if $u$ is a
bijection, and if the inequality $u(a)\leq u(b)$ is equivalent to the
inequality $a\leq b$ for any pair $a, b\in T_{D;W}$ (see [1, 3.5]).  We denote
by $Aut_0(T_{D;W})$ the group of automorphisms of the partially ordered set
$T_{D;W}$, such that $\alpha(T_{\mu;W})=T_{\mu;W}$ for any $\mu\in D$.  In
particular, when $D=P_d$, we obtain the group $Aut_0(T_{d;W})$.

Let $W\leq W^\prime\leq S_d$ be a pair of permutations groups such that
$|W^\prime:W|=2$.  In particular, $W$ is a normal subgroup of $W^\prime$ and
the factor-group $W^\prime /W$ has order $2$.  We denote by $\chi_e$ the
one-dimensional complex-valued character of $W^\prime$ with kernel $W$ (cf.
6.2).  In keeping with 6.2, for any $\lambda\in P_d$ the subset
$T_{\lambda;\chi_e}\subset T_{\lambda;W^\prime}$ consists of those
$W^\prime$-orbits in $T_\lambda$ which contain two $W$-orbits.  Let us set
$T_{D;\chi_e}=\cup_{\mu\in D}T_{\mu;\chi_e}$.  We fix a $\tau\in
W^\prime\backslash W$, so $\tau^2\in W$.  Then we have
$W^\prime/W=\langle\bar\tau\rangle$, where $\bar\tau=\tau W$.  The factor-group
$W^\prime /W$ acts on the set $T_{d;W}$ via the rule $(\eta W)O_W(A)=O_W(\eta
A)$.

We remind that for any $W\leq S_d$ the canonical projection $T_d\to T_{d;W}$ is
denoted by $\psi_W$ (see 4.1); its restriction on $T_D$ is the canonical
projection $\psi_{D;W}\colon T_D\to T_{D;W}$.  The following technical lemma is
obvious.

\proclaim{Lemma 8.1.1} Let $D\subset P_d$.

(i) The surjective map
$$
\psi_D\colon T_{D;W}\to
T_{D;W^\prime},\hbox{\ } O_W(A)\mapsto O_{W^\prime}(A),
$$
where $A\in T_D$, is the unique map such that
$\psi_D\psi_{D;W}=\psi_{D;W^\prime}$;

(ii) the canonical map
$$
(W^\prime/W)\backslash T_{D;W}\to T_{D;W^\prime},\hbox{\ }
O_{W^\prime/W}(O_W(A))\mapsto \cup_{\eta\in W^\prime}(\eta W)O_W(A)
$$
is a bijection;

(iii) after the identification of the orbit space $(W^\prime/W)\backslash
T_{D;W}$ with $T_{D;W^\prime}$ via the canonical bijection from (ii), the map
$\psi_D$ coincides with the canonical projection
$T_{D;W}\to (W^\prime/W)\backslash T_{D;W}$;

(iv) for $a^\prime\in T_{D;W^\prime}$ one has $\#\psi_D^{-1}(a^\prime)=2$ if
$a^\prime\in T_{D;\chi_e}$ and $\#\psi_D^{-1}(a^\prime)=1$ if
$a^\prime\in T_{D;W^\prime}\backslash T_{D;\chi_e}$;

(v) for any $W^\prime/W$-invariant map $\alpha\colon T_{D;W}\to T_{D;W}$ there
is a unique map $\alpha^\prime\colon T_{D;W^\prime}\to T_{D;W^\prime}$,
such that $\alpha^\prime\psi_D=\psi_D\alpha$;

(vi) for any two $W^\prime/W$-invariant maps $\alpha\colon T_{D;W}\to T_{D;W}$
and $\beta\colon T_{D;W}\to T_{D;W}$, one has
$(\alpha\beta)^\prime=\alpha^\prime\beta^\prime$.

\endproclaim

We denote by $Aut_0^{W^\prime}(T_{D;W})$ the subgroup of $Aut_0(T_{D;W})$
consisting of the $W^\prime/W$-invariant automorphisms $\alpha$, that is,
$\alpha(\iota a)=\iota\alpha(a)$ for any $\iota\in W^\prime/W$ and for any
$a\in T_{D;W}$.

The lemma below shows, that here we use correctly the notation
$Aut_0^{W^\prime}(T_{D;W})$ for the group of chiral automorphisms introduced in
Section 7.

\proclaim{Lemma 8.1.2} Let $\alpha\in Aut_0(T_{D;W})$.  The following three
statements are equivalent:

(i) one has $\alpha\in Aut_0^{W^\prime}(T_{D;W})$;

(ii) there exists an automorphism $\alpha^\prime\in Aut_0(T_{D;W^\prime})$,
such that $\psi_D\alpha=\alpha^\prime\psi_D$.

(iii) for any $a^\prime\in T_{D;\chi_e}$ with
$\psi^{-1}(a^\prime)=\{a,a_1\}$, $a\neq a_1$, there exists
$a^{\prime\prime}\in T_{D;\chi_e}$ such that
$\psi_D^{-1}(a^{\prime\prime})=\{\alpha(a),\alpha(a_1)\}$;

\endproclaim

\demo{Proof} (i) implies (ii).  In compliance with Lemma 8.1.1, (v), any
$W^\prime/W$-invariant map $\alpha\colon T_{D;W}\to T_{D;W}$ produces a unique
map $\alpha^\prime\colon T_{D;W^\prime}\to T_{D;W^\prime}$, such that
$\alpha^\prime\psi_D=\psi_D\alpha$.  The functoriality of this dependence
(Lemma 8.1.1, (vi)) implies that since $\alpha$ is a bijection, then
$\alpha^\prime$ is a bijection.  Now, we shall show that $\alpha^\prime\in
Aut_0(T_{D;W^\prime})$.  Obviously,
$\alpha^\prime(T_{\mu;W^\prime})=T_{\mu;W^\prime}$ for all $\mu\in P_d$.
Suppose that $a^\prime\leq \bar a^\prime$ in $T_{D;W^\prime}$.  Equivalently,
there are $a,\bar a\in T_{D;W}$ such that $a\subset a^\prime$, $\bar a\subset
\bar a^\prime$, and $a\leq\bar a$.  The last inequality s equivalent to
$\alpha(a)\leq\alpha(\bar a)$, and since
$\alpha(a)\subset\alpha^\prime(a^\prime)$, $\alpha(\bar
a)\subset\alpha^\prime(\bar a^\prime)$, we obtain
$\alpha^\prime(a^\prime)\leq\alpha^\prime(\bar a^\prime)$.  Conversely, the
definition of $\alpha^\prime$ and the previous inequality imply the existence
of $a,\bar a\in T_{D;W}$ such that $\alpha(a)\leq\alpha(\bar a)$.  Hence the
two inequalities $\alpha^\prime(a^\prime)\leq\alpha^\prime(\bar a^\prime)$ and
$a^\prime\leq\bar a^\prime$ are equivalent.

(ii) implies (iii).  Suppose that $a^\prime\in T_{D;\chi_e}$. Then we have
$\psi_D^{-1}(a^\prime)=\{a,a_1\}$, where $a\neq a_1$.  We set
$a^{\prime\prime}= \alpha^\prime(a^\prime)$.  Thus, $a^{\prime\prime}=
\alpha^\prime(\psi_D(a))=\psi_D(\alpha(a))$, and, by analogy,
$a^{\prime\prime}=\psi_D(\alpha(a_1))$, so
$\psi_D^{-1}(a^{\prime\prime})=\{\alpha(a),\alpha(a_1)\}$.  In particular,
$a^{\prime\prime}\in T_{D;\chi_e}$.

(iii) implies (i).  Suppose that (iii) holds.  Then $\alpha$ maps the set
$\psi_D^{-1}(T_{D;\chi_e})$ and its complement onto themselves.  We want to
show that $\alpha(\iota a)=\iota\alpha(a)$ for any $\iota\in W^\prime/W$ and
for any $a\in T_{D;W}$.  If $\iota$ is the unit element of $W^\prime/W$, our
equality is trivial.  Now, let $\iota=\bar\tau$.  If $a\neq\bar\tau a$, then
$\psi_D^{-1}(a^\prime)=\{a,\bar\tau a\}$, where $a^\prime=\psi_D(a)$.
Therefore $\psi_D^{-1}(a^{\prime\prime})=\{\alpha(a),\alpha(\bar\tau a)\}$, and
hence $\alpha(\bar\tau a)=\bar\tau\alpha(a)$.  If $a=\bar\tau a$, then
$\bar\tau\alpha(a)=\alpha(a)$, and we have again $\alpha(\bar\tau
a)=\bar\tau\alpha(a)$.

\enddemo

8.2. We remind that via the rule (1.1.1) any permutation $\zeta\in S_d$ induces
a bijection $\zeta\colon T_d\to T_d$, $A\mapsto \zeta A$, of the set of all
tabloids with $d$ nodes, and its restriction on $T_\mu$ is a bijection of
$T_\mu$ onto itself for all $\mu\in P_d$.

Let $W\leq W^\prime\leq S_d$ be a pair of permutation groups such that
$|W^\prime:W|=2$.  We denote by $N$ the normalizer of the group $W$ in $S_d$,
and by $N^\prime$ the intersection of $N$ with the normalizer of $W^\prime$ in
$S_d$.  Plainly, $W^\prime\leq N^\prime$.

\proclaim{Lemma 8.2.1} For any $\nu\in N$ the bijection $\nu$ can be factored
out to a bijection $\hat{\nu}\colon T_{d;W}\to T_{d;W}$, $a\mapsto
\hat{\nu}(a)$:  if $A\in a$, then $\hat{\nu}(a)$ is the $W$-orbit of $\nu A$.
Moreover, the bijection $\hat{\nu}$ maps the set $T_{\lambda;W}$ onto itself
for any $\lambda\in P_d$.

\endproclaim

\demo{Proof} Given $\nu\in N$ and $\sigma\in W$, there exists $\sigma_1\in W$,
such that $\nu\sigma=\sigma_1\nu$.  Therefore $O_W(\nu\sigma A)=
O_W(\sigma_1\nu A)= O_W(\nu A)$, so the rule in the condition yields a map
$\hat\nu$.  We have $(\hat\nu)^{-1}=\widehat{\nu^{-1}}$, hence $\hat\nu$ is a
bijection of $T_{d;W}$ onto itself, and plainly any set $T_{\lambda;W}$ is
$\hat\nu$-stable.

\enddemo

\proclaim{Theorem 8.2.2} (i) Given $\nu\in N$ and $D\subset P_d$, the bijection
$\hat{\nu}$ induces an automorphism from $Aut_0(T_{D;W})$;

(ii) the map $\nu\mapsto\hat{\nu}$ is a homomorphism of the group $N$ into the
group $Aut_0(T_{D;W})$, whose kernel contains $W$; if $(1^d)\in D$, then its
kernel coincides with $W$;

(iii) the map $\nu\mapsto\hat{\nu}$ is a homomorphism of the group $N^\prime$
into the group $Aut_0^{W^\prime}(T_{D;W})$, whose kernel contains $W$; if
$(1^d)\in D$, then its kernel coincides with $W$.

\endproclaim

\demo{Proof} (i) According to Lemma 8.2.1, $\hat{\nu}$ induces a bijection of
the set $T_{D;W}$ onto itself.  Let $a,b\in T_{D;W}$, and $A\in a$, $B\in b$.
The inequality $a\leq b$ means that there exists a $\sigma\in W$, such that
$\sigma A\leq B$, or, equivalently, $\nu\sigma A\leq\nu B.$ The inequality
$\hat{\nu}(a)\leq\hat{\nu}(b)$ means that there exists a $\sigma_1\in W$, such
that $\sigma_1\nu A\leq\nu B.$ Since $\nu W=W\nu$, the last two conditions are
equivalent.

(ii) Let $\nu_1,\nu_2\in N$, and let $a\in T_{D;W}$, $a=O_W(A)$. We have
$$
\widehat{\nu_1\nu_2}(a)=
O_W(\nu_1\nu_2A)=
O_W(\nu_1(\nu_2A))=
\hat{\nu_1}(O_W(\nu_2A))=
\hat{\nu_1}(\hat{\nu_2}(a))=
(\hat{\nu_1}\hat{\nu_2})(a),
$$
so
$\widehat{\nu_1\nu_2}=\hat{\nu_1}\hat{\nu_2}$.

Plainly, if $\sigma\in W$, then $\hat{\sigma}(a)=a$ for any $a\in T_{D;W}$,
hence $W$ is a subgroup of the kernel.  Now, suppose $(1^d)\in D$, and let
$\nu\in N$ be such that $\hat{\nu}(a)=a$ for any $a\in T_{D;W}$.  The last
condition can also be written as $O_W(\nu A)=O_W(A)$ for any $A\in T_D$.
In other words, for any $A=(A_1,A_2,\hdots)\in
T_D$ there exists a $\sigma\in W$ such that $\nu A=\sigma A$, that is,
$\nu(A_1)=\sigma(A_1)$, $\nu(A_2)=\sigma(A_2),\hdots$.  In particular, if
$A=(A_1,A_2,\hdots, A_d)$ with $A_1=\{1\}, A_2=\{2\},\hdots ,A_d=\{d\}$, then
$\nu=\sigma\in W$.

(iii) Part (i) applied for the group $W^\prime$ yields that any $\nu\in
N^\prime$ induces an automorphism ${\hat\nu}^\prime\in Aut_0(T_{D;W^\prime})$.
It is obvious that $\psi_D\hat\nu={\hat\nu}^\prime\psi_D$. Then Lemma 8.1.2
implies
$\hat\nu\in Aut_0^{W^\prime}(T_{D;W})$. Now, since $N^\prime$ is a subgroup
of $N$ that contains $W$, part (ii) finishes the proof.

\enddemo

Theorem 8.2.2, (ii), yields information concerning the subgroup
$Aut_0^{h;W^\prime}(T_{D;W})$ of the group $Aut_0^{W^\prime}(T_{D;W})$,
consisting of all hidden symmetries.

\proclaim{Corollary 8.2.3} If $(1^d)\in D$, then there exists a natural
embedding of the factor-group $N^\prime/W$ into the group
$Aut_0^{W^\prime}(T_{D;W})$ with image $Aut_0^{h;W^\prime}(T_{D;W})$.

\endproclaim

We set $D_e=\varphi_W^\prime(T_{d;\chi_e})$, where $\varphi_W^\prime\colon
T_{d;W^\prime}\to P_d$ is the restriction of the projection defined in the
beginning of Subsection 4.1.  Theorem 5.3.1 yields immediately that the subset
$D_e\subset P_d$ is an order ideal (see [1, Section 2]) of the partially
ordered set $P_d$.  Moreover, $D_e\neq\emptyset$ because at least
$T_{\left(1^d\right);W^\prime}\subset T_{d;\chi_e}$.

\proclaim{Corollary 8.2.4} For any $D\subset P_d$ there exists a natural
homomorphism of the factor-group $W^\prime/W=\langle\bar\tau\rangle$ into the
group $Aut_0^{W^\prime}(T_{D;W})$ with image $\langle\hat\tau\rangle$.  Given
$a^\prime\in T_{D;W^\prime}$, the automorphism $\hat\tau \in
Aut_0^{W^\prime}(T_{D;W})$ maps the fiber $\psi^{-1}(a^\prime)$ onto itself in
the following way:  if $a^\prime\in T_{D;W^\prime}\backslash T_{D;\chi_e}$,
then $\psi^{-1}(a^\prime)=\{a\}$ and $\hat\tau(a)=a$, and if $a^\prime\in
T_{D;\chi_e}$, then $\psi^{-1}(a^\prime)=\{a,a_1\}$ and $\hat\tau(a)=a_1$,
$\hat\tau(a_1)=a$.  If $D\cap D_e=\emptyset$, then $\hat\tau$ is the identity.
If $D\cap D_e\neq\emptyset$, then $\hat\tau$ is an involution and the above
homomorphism is an embedding $W^\prime/W\simeq\langle\hat\tau\rangle\leq
Aut_0^{h;W^\prime}(T_{D;W})$.

\endproclaim

\demo{Proof} Because of $W^\prime\leq N^\prime$, the homomorphism from Theorem
8.2.2, (iii), induces a homomorphism $W^\prime/W\to Aut_0^{W^\prime}(T_{D;W})$,
$\bar\eta\mapsto \hat\eta$, where $\bar\eta=\eta W$, $\eta\in W^\prime$.  In
particular, since $\bar\eta a=\hat\eta (a)$ for any $a\in T_{d;W}$, $\eta\in
W^\prime$, Lemma 8.1.1, (iv), yields the prescribed action on the fibers of
$\psi$.  If $D\cap D_e=\emptyset$, then for all $a^\prime\in T_{D;W^\prime}$
the fiber $\psi^{-1}(a^\prime)$ consist of one inverse image $a$:
$\psi^{-1}(a^\prime)=\{a\}$, and this yields $\hat\tau(a)=a$, so $\hat\tau$ is
the identity automorphism.  The condition $D\cap D_e\neq\emptyset$ implies that
there exists at least one $a^\prime\in T_{D;W^\prime}$ with two inverse images:
$\psi^{-1}(a^\prime)=\{a,a_1\}$.  Then $\hat\tau(a)=a_1$, $\hat\tau(a_1)=a$,
and as a consequence the automorphism $\hat\tau$ is an involution.  Therefore
we get the desired embedding.

\enddemo

In accord with Lemma 8.1.1, (v), and Lemma 8.1.2, for any $\alpha\in
Aut_0^{W^\prime}(T_{d;W})$ there exists a unique automorphism $\alpha^\prime\in
Aut_0(T_{d;W^\prime})$, such that $\psi\alpha=\alpha^\prime\psi$.

\proclaim{Proposition 8.2.5} For any $D\subset P_d$ the map
$$
Aut_0^{W^\prime}(T_{D;W})\to
Aut_0(T_{D;W^\prime}),\hbox{\ }\alpha\mapsto\alpha^\prime,
$$
is a homomorphism of groups and the group $\langle\hat\tau\rangle$ is contained
in its kernel.

\endproclaim

\demo{Proof} Let $\alpha,\beta\in Aut_0^{W^\prime}(T_{D;W})$.  Lemma 8.1.1,
(vi), yields $(\alpha\beta)^\prime=\alpha^\prime\beta^\prime$.  Moreover,
$(\hat\tau)^\prime$ is the identity.

\enddemo

The elements of the group $Aut_0^{W^\prime}(T_{D;W})$ are called {\it chiral
automorphisms} of the partially ordered set $T_{D;W}$.  In case $D\cap
D_e\neq\emptyset$ the involution $\hat\tau\in Aut_0^{W^\prime}(T_{D;W})$ is
said to be the {\it chiral involution} of the group
$Aut_0^{W^\prime}(T_{D;W})$.

When $W=G$ and $W^\prime=G^\prime$, where $G$ and $G^\prime$ are the groups
of substitution isomerism and stereoisomerism, respectively, of the molecule
under consideration, and when a distribution of ligants amounts to a chiral
molecule, then in compliance with (1.5.1), (2a), we have $|G^\prime :G|=2$ and
denote the group $Aut_0^{G^\prime}(T_{D;G})$ by $Aut_0^\prime(T_{D;G})$.

Taking into account Lunn-Senior's thesis 1.5.1, (2b), we can reformulate
(8.2.4) in the following way:

\proclaim{Theorem 8.2.6} The members of a chiral pair are indistinguishable via
substitution reactions.

\endproclaim

\heading
9. Characters

\endheading

9.1. Here we continue to use notations, introduced in Section 8.
Let $X_W$ be the Abelian group of the one-dimensional characters of the
group $W$. The normalizer $N$ of $W$ in $S_d$ acts on $X_W$ by the rule
$$
(\nu\chi)(\sigma)=\chi(\nu^{-1}\sigma\nu),
$$
where $\chi\in X_W$, $\sigma\in W$, $\nu\in N$.
Since $\sigma\chi=\chi$ for each $\sigma\in W$, we obtain an action of the
factor-group $N/W$ on $X_W$, defined by the rule
$$
(\nu W\chi)(\sigma)=\chi(\nu^{-1}\sigma\nu),
$$
where $\chi\in X_W$, $\nu\in N$. In particular, since $N^\prime/W$ is a
subgroup of $N/W$, and because of Corollary 8.2.3, in case $(1^d)\in D$ we
obtain an action of the group $Aut_0^{h;W^\prime}(T_{D;W})$ of hidden
symmetries on the set $X_W$:
$$
(\alpha\chi)(\sigma)=\chi(\nu^{-1}\sigma\nu),
$$
where $\chi\in X_W$, $\alpha\in Aut_0^{h;W^\prime}(T_{D;W})$, $\alpha=\hat\nu$,
$\nu\in N^\prime$.

Now, we refer to the notations from Subsection 5.1.  Let $\lambda$ be a
partition of the number $d$ and let $\theta$ be a one-dimensional character of
the group $S_\lambda$.  Let $\chi\in X_W$.  Given a tabloid $A=\upsilon I\in
T_\lambda$, we note that the one-dimensional character $\beta_\upsilon$ of the
stabilizer $W_A$ from (5.1.2) depends only on $A$, and on $\chi$ when $\theta$
is fixed, so we introduce the more precise notation $\beta_\upsilon=
\beta_{A,\chi}$.

\proclaim{Theorem 9.1.1} Given $\nu\in N$, $\lambda\in P_d$, $\chi\in X_W$,
$\theta\in X_{S_\lambda}$, one has:

(i) the automorphism $\hat{\nu}$ of the partially ordered set $T_{d;W}$ maps
the set $T_{\lambda;\chi,\theta}$ of all $(\chi,\theta)$-orbits onto the set
$T_{\lambda;\nu\chi,\theta}$ of all $(\nu\chi,\theta)$-orbits;

(ii) if $(1^d)\in D$, the hidden automorphism $\alpha\in
Aut_0^{h;W^\prime}(T_{D;W})$ of the
partially ordered set $T_{D;W}$ maps the set $T_{\lambda;\chi,\theta}$ of all
$(\chi,\theta)$-orbits in $T_{D;W}$ onto the set
$T_{\lambda;\alpha\chi,\theta}$ of all $(\alpha\chi,\theta)$-orbits in
$T_{D;W}$.

\endproclaim

\demo{Proof} Obviously part (ii) is a consequence of part (i).

(i) Let $a\in T_{\lambda;\chi,\theta}$, and let $A\in a$, $A=\upsilon I$.  For
any $\nu\in N$ we have $\nu A=\nu\upsilon I$, $W_{\nu A}=\nu W_A\nu^{-1}$, and
if $\sigma\in W$, then
$$
\beta_{\nu A,\nu\chi}(\nu\sigma\nu^{-1})=
(\nu\chi)(\nu\sigma\nu^{-1})\theta((\nu\upsilon)^{-1}\nu\sigma\nu^{-1}\nu
\upsilon)=
\chi(\sigma)\theta(\upsilon^{-1}\sigma\upsilon)=\beta_{A,\chi}(\sigma).
$$
In particular, $\beta_{A,\chi}=1$ on $W_A$ if and only if $\beta_{\nu
A,\nu\chi}=1$ on $W_{\nu A}$.  Therefore if $a\in T_{\lambda;\chi,\theta}$,
then $\hat{\nu}(a)\in T_{\lambda;\nu\chi,\theta}$, and if $a\in
T_{\lambda;\nu\chi,\theta}$, then $\hat{\nu}^{-1}(a)= \widehat{\nu^{-1}}(a)\in
T_{\lambda;\chi,\theta}$, and we are done.

\enddemo

By substituting $\theta=1_{S_\lambda}$ in (9.1.1) we obtain

\proclaim{Corollary 9.1.2} For any $\nu\in N$, any $\lambda\in
S_\lambda$, and any $\chi\in X_W$ one has:

(i) the automorphism $\hat{\nu}\in Aut_0(T_{d;W})$ maps the set
$T_{\lambda;\chi}$ of all $\chi$-orbits onto the set $T_{\lambda;\nu\chi}$ of
all $\nu\chi$-orbits;

(ii) if $(1^d)\in D$, the hidden automorphism $\alpha\in
Aut_0^{h;W^\prime}(T_{D;W})$ maps the set $T_{\lambda;\chi}$ of all
$\chi$-orbits in $T_{D;W}$ onto the set $T_{\lambda;\alpha\chi}$ of all
$\alpha\chi$-orbits in $T_{D;W}$.

\endproclaim

\proclaim{Corollary 9.1.3} For any $\nu\in N$, $\lambda\in P_d$, $\chi\in X_W$,
and $\theta\in X_{S_\lambda}$, one has:

(i) $n_{\lambda;\nu\chi,\theta}=
n_{\lambda;\chi,\theta}$ and $n_{\lambda;\nu\chi}= n_{\lambda;\chi}$;

(ii) if $(1^d)\in D$, then for any hidden automorphism $\alpha\in
Aut_0^{h;W^\prime}(T_{D;W})$ one has $n_{\lambda;\alpha\chi,\theta}=
n_{\lambda;\chi,\theta}$ and $n_{\lambda;\alpha\chi}= n_{\lambda;\chi}$.

\endproclaim

Taking into account (9.1.3), (i), (5.2.4), (i), (iv), and [8, Ch. I, Sec. 7,
7.3], we establish

\proclaim{Corollary 9.1.4} (i) For any $\nu\in N$, and any $\chi\in X_W$ the
induced monomial representations $Ind_W^{S_d}(\chi)$ and $Ind_W^{S_d}(\nu\chi)$
of the symmetric group $S_d$, are equivalent;

(ii) if $(1^d)\in D$, then for any hidden automorphism $\alpha\in
Aut_0^{h;W^\prime}(T_{D;W})$ the induced monomial representations
$Ind_W^{S_d}(\chi)$ and $Ind_W^{S_d}(\alpha\chi)$ of the symmetric group $S_d$,
are equivalent.

\endproclaim

9.2.  Let $a, b\in T_{\lambda;W}$, $\lambda\in P_d$, and let $\chi\in X_W$,
$\theta\in X_{S_{\lambda}}$.  It is said that the pair of characters
$(\chi,\theta)$ {\it separates $a$ from $b$} if $a\in T_{\lambda;\chi,\theta}$
and $b\notin T_{\lambda;\chi,\theta}$.  The character $\chi\in X_W$ {\it
separates $a$ from $b$} if $a\in T_{\lambda;\chi}$ and $b\notin
T_{\lambda;\chi}$.  We say that $a$ and $b$ are
{\it indistinguishable via pairs of characters}
if:  (a$_{PC}$) $a$ and $b$ are contained in one and the same
$G^{\prime\prime}$-orbit, and (b$_{PC}$) for any pair of characters
$(\chi,\theta)$ that separates $a$ from $b$ there exists a permutation $\nu\in
N^\prime$ such that $(\nu\chi,\theta)$ separates $b$ from $a$, and for any pair
of characters $(\chi,\theta)$ that separates $b$ from $a$ there exists a
permutation $\nu\in N^\prime$ such that $(\nu\chi,\theta)$ separates $a$ from
$b$.

It is said that $a$ and $b$ are {\it indistinguishable via characters} if:
(a$_C$) $a$ and $b$ are contained in one and the same $G^{\prime\prime}$-orbit,
and (b$_C$) for any character $\chi\in X_W$ that separates $a$ from $b$ there
exists a permutation $\nu\in N^\prime$ such that $\nu\chi$ separates $b$ from
$a$, and for any character $\chi\in X_W$ that separates $b$ from $a$ there
exists a permutation $\nu\in N^\prime$ such that $\nu\chi$ separates $a$ from
$b$.

Otherwise, $a$ and $b$ are said to be {\it distinguishable via pairs of
characters}, respectively, {\it via characters}.

By setting $\theta=1_{S_\lambda}$ we obtain that if $a$ and $b$ are
indistinguishable via pairs of characters, then $a$ and $b$ are indistinguishable
via characters.

\proclaim{Theorem 9.2.1} The members of a chiral pair are indistinguishable
via pairs of characters.

\endproclaim

\demo{Proof} Let $a^\prime\in T_{\lambda;W^\prime}$ be a chiral pair, that is,
$\psi^{-1}(a^\prime)=\{a,a_1\}$, $a,a_1\in T_{\lambda;W}$, $a\neq a_1$.  In
compliance with Lunn-Senior's thesis 1.5.1, (2b), $a$ and $a_1$ are contained
in one and the same $G^{\prime\prime}$-orbit.  Corollary 8.2.4 yields that the
chiral involution $\hat\tau\in Aut_0^{h;W^\prime}(T_{d;W})$ permutes $a$ and
$a_1$.  Now, suppose that the pair of characters $(\chi,\theta)$ separates
$a$ from $a_1$, so $a\in T_{\lambda;\chi,\theta}$ and $a_1\notin
T_{\lambda;\chi,\theta}$.  Then in compliance with Theorem 9.1.1, (i),
$\hat\tau(T_{\lambda;\chi,\theta})=T_{\lambda;\tau\chi,\theta}$.  Therefore
$a_1\in T_{\lambda;\tau\chi,\theta}$ and $a\notin T_{\lambda;\tau\chi,\theta}$,
so the pair $((\tau\chi,\theta)$ separates $a_1$ from $a$.  Thus, $a$ and $a_1$
are not distinguishable via pairs of characters.

\enddemo

\proclaim{Corollary 9.2.2} The members of a chiral pair are indistinguishable
via characters.

\endproclaim

\heading
10. Examples of groups of chiral automorphisms \\

\endheading

10.1.  In this Subsection our constant reference is 6.4.  We remind that the
group $G$ of univalent substitution isomerism of ethene $C_2H_4$ is the Klein
four group $G=\{(1), (12)(34), (13)(24), (14)(23)\}$ in $S_4$.  There are no
chiral pairs among the derivatives of ethene, so $G=G^\prime$, and hence the
group $Aut_0^\prime(T_{4;G})$ of chiral automorphisms coincides with the group
$Aut_0(T_{4;G})$.

Here is the complete diagram:

\vskip 12 pt

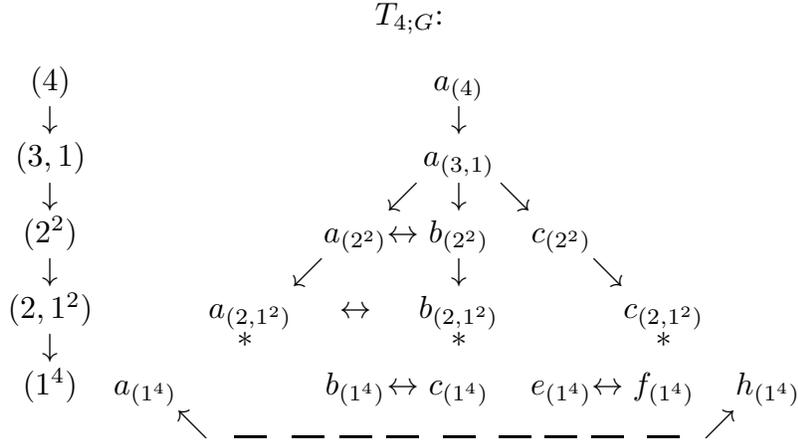
\centerline{ $T_{4;G}$:}

$$
\left.
\matrix \format \c & \c & \c & \c & \c & \c & \c & \c & \c & \c & \c & \c & \c
& \c &  \c  \\
(4) & \hbox{\ \ } & \hbox{\ \ } & \hbox{\ \ } & \hbox{\ \ } & \hbox{\ \ } &
\hbox{\ \ } & \hbox{\ \ } & a_{\left(4\right)} & \hbox{\ \ } & \hbox{\ \ } &
\hbox{\ \ } & \hbox{\ \ } & \hbox{\ \ } & \hbox{\ \ } \cr
\downarrow & \hbox{\ \ } & \hbox{\ \ } & \hbox{\ \ } & \hbox{\ \ } & \hbox{\ \
} & \hbox{\ \ } & \hbox{\ \ } & \downarrow & \hbox{\ \ } & \hbox{\ \ } &
\hbox{\ \ } & \hbox{\ \ } & \hbox{\ \ } & \hbox{\ \ } \cr
(3,1) & \hbox{\ \ } & \hbox{\ \ } & \hbox{\ \ } & \hbox{\ \ } & \hbox{\ \ } &
\hbox{\ \ } & \hbox{\ \ } & a_{\left(3,1\right)} & \hbox{\ \ } & \hbox{\ \ } &
\hbox{\ \ } & \hbox{\ \ } & \hbox{\ \ } & \hbox{\ \ } \cr
\downarrow & \hbox{\ \ } & \hbox{\ \ } & \hbox{\ \ } & \hbox{\ \ } & \hbox{\ \
} & \hbox{\ \ } & \swarrow & \downarrow & \searrow & \hbox{\ \ } & \hbox{\ \ }
& \hbox{\ \ } & \hbox{\ \ } & \hbox{\ \ } \cr
(2^2) & \hbox{\ \ } & \hbox{\ \ } & \hbox{\ \ } & \hbox{\ \ } & \hbox{\ \ } &
a_{\left(2^2\right)} & \leftrightarrow &
b_{\left(2^2\right)} & \hbox{\ \ } & c_{\left(2^2\right)} & \hbox{\ \ } &
\hbox{\ \ } & \hbox{\ \ } & \hbox{\ \ } \cr
\downarrow & \hbox{\ \ } & \hbox{\ \ } & \hbox{\ \ } & \hbox{\ \ } & \swarrow
& \hbox{\ \ } & \hbox{\ \ } & \downarrow & \hbox{\ \ } & \hbox{\ \ } & \searrow
& \hbox{\ \ } & \hbox{\ \ } & \hbox{\ \ } \cr
(2,1^2) & \hbox{\ \ } & \hbox{\ \ } & \hbox{\ \ } & a_{\left(2,1^2\right)} &
\hbox{\ \ } & \leftrightarrow &
\hbox{\ \ } & b_{\left(2,1^2\right)} & \hbox{\ \ } & \hbox{\ \ } & \hbox{\ \ }
& c_{\left(2,1^2\right)} & \hbox{\ \ } & \hbox{\ \ } \cr
\downarrow & \hbox{\ \ } & \hbox{\ \ } & \hbox{\ \ } & \hbox{* } & \hbox{\ \ }
& \hbox{\ \ } & \hbox{\ \ } & \hbox{*} & \hbox{\ \ } & \hbox{\ \ } & \hbox{\ \ }
& \hbox{*} & \hbox{\ \ } & \hbox{\ \ } \cr
(1^4) & \hbox{\ \ } & a_{\left(1^4\right)} & \hbox{\ \ } & \hbox{\ \ } &
\hbox{\ \ } & b_{\left(1^4\right)} &
\leftrightarrow & c_{\left(1^4\right)} & \hbox{\ \ } & e_{\left(1^4\right)}
&\leftrightarrow & f_{\left(1^4\right)}
& \hbox{\ \ } & h_{\left(1^4\right)} \cr
\hbox{\ \ } & \hbox{\ \ } & \hbox{\ \ } & \nwarrow & \underline{\ \ \ } &
\underline{\ \ \ } & \underline{\ \ \ } & \underline{\ \ \ } & \underline{\ \ \
} & \underline{\ \ \ } & \underline{\ \ \ } & \underline{\ \ \ } & \underline{\
\ \ } & \nearrow & \hbox{\ \ } \cr
\endmatrix \right.
$$

\vskip 12 pt

The horizontal double arrow means that the two diamers are identical as
structural isomers.  The star $*$ below a structural formula indicates that all
possible arrows from it to any structural formula from the next row are drawn.

Plainly, any automorphism $\alpha$ from $Aut_0(T_{4;G})$ leaves invariant the
elements $a_{\left(4\right)}$ and $a_{\left(3,1\right)}$, and if, for instance,
$\alpha(a_{\left(2^2\right)})=b_{\left(2^2\right)}$, then
$\alpha(a_{\left(2,1^2\right)})=b_{\left(2,1^2\right)}$.  Similarly, if for
example, $\alpha(a_{\left(2^2\right)})=a_{\left(2^2\right)}$, then
$\alpha(a_{\left(2,1^2\right)})=a_{\left(2,1^2\right)}$, etc.  On the six
letters $a_{\left(1^4\right)}$, $b_{\left(1^4\right)}$, $c_{\left(1^4\right)}$,
$e_{\left(1^4\right)}$, $f_{\left(1^4\right)}$, $h_{\left(1^4\right)}$, there
are no restrictions on $\alpha$.  Thus, we have proved

\proclaim{Theorem 10.1.1} One has
$$
Aut_0(T_{4;G})\simeq S_3\times S_6,
$$
where the permutations of $S_3$ simultaneously move the elements $a$, $b$, $c$
with indices $(2^2)$ and $(2,1^2)$, and $S_6$ consists of all permutations of
the six elements with indices $(1^4)$.

\endproclaim

We have $N=N^\prime=S_4$, so in accord to Corollary 8.2.3 there are
$|N/G|=6$ hidden symmetries in $Aut_0(T_{4;G})$, and we obtain

\proclaim{Proposition 10.1.2} The group $Aut_0^h(T_{4;G})$ of hidden symmetries
of ethene consists of the unit permutation, and
$$
\widehat{(12)}=
(b_{\left(2^2\right)}\hbox{\ }c_{\left(2^2\right)})
(b_{\left(2,1^2\right)}\hbox{\ }c_{\left(2,1^2\right)})
(a_{\left(1^4\right)}\hbox{\ }b_{\left(1^4\right)})
(c_{\left(1^4\right)}\hbox{\ }e_{\left(1^4\right)})
(f_{\left(1^4\right)}\hbox{\ }h_{\left(1^4\right)}),
$$
$$
\widehat{(23)}=
(a_{\left(2^2\right)}\hbox{\ }c_{\left(2^2\right)})
(a_{\left(2,1^2\right)}\hbox{\ }c_{\left(2,1^2\right)})
(a_{\left(1^4\right)}\hbox{\ }e_{\left(1^4\right)})
(b_{\left(1^4\right)}\hbox{\ }f_{\left(1^4\right)})
(c_{\left(1^4\right)}\hbox{\ }h_{\left(1^4\right)}),
$$
$$
\widehat{(13)}=
(a_{\left(2^2\right)}\hbox{\ }b_{\left(2^2\right)})
(a_{\left(2,1^2\right)}\hbox{\ }b_{\left(2,1^2\right)})
(a_{\left(1^4\right)}\hbox{\ }h_{\left(1^4\right)})
(b_{\left(1^4\right)}\hbox{\ }c_{\left(1^4\right)})
(e_{\left(1^4\right)}\hbox{\ }f_{\left(1^4\right)}),
$$
$$
\widehat{(123)}=
(a_{\left(2^2\right)}\hbox{\ }b_{\left(2^2\right)}
\hbox{\ }c_{\left(2^2\right)})
(a_{\left(2,1^2\right)}\hbox{\ }b_{\left(2,1^2\right)}
\hbox{\ }c_{\left(2,1^2\right)})
(a_{\left(1^4\right)}\hbox{\ }c_{\left(1^4\right)}
\hbox{\ }f_{\left(1^4\right)})
(b_{\left(1^4\right)}\hbox{\ }h_{\left(1^4\right)}
\hbox{\ }e_{\left(1^4\right)}),
$$
$$
\widehat{(124)}=
(a_{\left(2^2\right)}\hbox{\ }c_{\left(2^2\right)}
\hbox{\ }b_{\left(2^2\right)})
(a_{\left(2,1^2\right)}\hbox{\ }c_{\left(2,1^2\right)}
\hbox{\ }b_{\left(2,1^2\right)})
(a_{\left(1^4\right)}\hbox{\ }f_{\left(1^4\right)}
\hbox{\ }c_{\left(1^4\right)})
(b_{\left(1^4\right)}\hbox{\ }e_{\left(1^4\right)}
\hbox{\ }h_{\left(1^4\right)}).
$$

\endproclaim

This yields

\proclaim{Corollary 10.1.3} The diamers of ethene in the pairs
$\{a_{\left(2^2\right)},b_{\left(2^2\right)}\}$,
$\{a_{\left(2,1^2\right)},b_{\left(2,1^2\right)}\}$, $\{a_{\left(1^4\right)},
h_{\left(1^4\right)}\}$, $\{b_{\left(1^4\right)},c_{\left(1^4\right)}\}$,
$\{e_{\left(1^4\right)}, f_{\left(1^4\right)}\}$
are indistinguishable via
substitution reactions.

\endproclaim

\demo{Proof} Proposition 10.1.2 implies that the permutation
$$
\widehat{(13)}=
(a_{\left(2^2\right)},b_{\left(2^2\right)})
(a_{\left(2,1^2\right)},b_{\left(2,1^2\right)})
(a_{\left(1^4\right)}, h_{\left(1^4\right)})
(b_{\left(1^4\right)},c_{\left(1^4\right)})(e_{\left(1^4\right)},
f_{\left(1^4\right)})
$$
is an element of the group $Aut_0^h(T_{4;G})\leq Aut_0^\prime(T_{4;G})$, and,
moreover, in agreement with the consideration in [7, VI] and in 6.4, the
diamers in the above pairs are structurally identical.

\enddemo

The Abelian group $X_G$ of one-dimensional characters of the group $G$ consists
of the unit character, and the characters $\chi_1$, $\chi_2$, and $\chi_3$,
defined in 6.4.  The normalizer $N=S_4$ acts on $X_G$ and the unit character is
one of the orbits of this action.  The other orbit is
$$
\{\chi_1,\chi_2,\chi_3\},
$$
because $\widehat{(123)}\chi_1=\chi_2$, and $\widehat{(123)}\chi_2=\chi_3$.

On the other hand, the Abelian group $X_{S_{\left(2^2\right)}}$ of
one-dimensional characters of the group
$S_{\left(2^2\right)}=\langle (12),(34)\rangle$ has four elements: the unit
character $\theta_0$, and the characters
$\theta_1$,
$\theta_2$, and
$\theta_3$,
defined by the equalities
$$
\theta_1((12))=1,\hbox{\ }
\theta_1((34))=-1,
$$
$$
\theta_2((12))=-1,\hbox{\ }
\theta_2((34))=1,
$$
and
$$
\theta_3((12))=-1,\hbox{\ }
\theta_3((34))=-1.
$$

The propositions below, together with (10.1.3), confirms the conclusions in [7,
VI], and in 6.4.

\proclaim{Proposition 10.1.4} The members of the diameric pairs
$\{a_{\left(2^2\right)},b_{\left(2^2\right)}\}$,
$\{a_{\left(2,1^2\right)},b_{\left(2,1^2\right)}\}$, $\{a_{\left(1^4\right)},
h_{\left(1^4\right)}\}$, $\{b_{\left(1^4\right)},c_{\left(1^4\right)}\}$,
$\{e_{\left(1^4\right)}, f_{\left(1^4\right)}\}$
cannot be distinguished via pairs of characters.

\endproclaim

\demo{Proof} In compliance with 5.1 and 6.4, we obtain

$$
a_{\left(2^2\right)}\notin T_{\left(2^2\right);\chi_1,\theta_0},\hbox{\ }
b_{\left(2^2\right)}\notin T_{\left(2^2\right);\chi_1,\theta_0},\hbox{\ }
a_{\left(2^2\right)}\in T_{\left(2^2\right);\chi_2,\theta_0},\hbox{\ }
b_{\left(2^2\right)}\notin T_{\left(2^2\right);\chi_2,\theta_0},
$$
$$
a_{\left(2^2\right)}\notin T_{\left(2^2\right);\chi_3,\theta_0}, \hbox{\ }
b_{\left(2^2\right)}\in T_{\left(2^2\right);\chi_3,\theta_0},
$$

$$
a_{\left(2^2\right)}\in T_{\left(2^2\right);\chi_1,\theta_1},\hbox{\ }
b_{\left(2^2\right)}\in T_{\left(2^2\right);\chi_1,\theta_1},\hbox{\ }
a_{\left(2^2\right)}\notin T_{\left(2^2\right);\chi_2,\theta_1},\hbox{\ }
b_{\left(2^2\right)}\in T_{\left(2^2\right);\chi_2,\theta_1},\hbox{\ }
$$
$$
a_{\left(2^2\right)}\in T_{\left(2^2\right);\chi_3,\theta_1},\hbox{\ }
b_{\left(2^2\right)}\notin T_{\left(2^2\right);\chi_3,\theta_1},
$$

$$
a_{\left(2^2\right)}\in T_{\left(2^2\right);\chi_1,\theta_2},\hbox{\ }
b_{\left(2^2\right)}\in T_{\left(2^2\right);\chi_1,\theta_2},\hbox{\ }
a_{\left(2^2\right)}\notin T_{\left(2^2\right);\chi_2,\theta_2},\hbox{\ }
b_{\left(2^2\right)}\in T_{\left(2^2\right);\chi_2,\theta_2},
$$
$$
a_{\left(2^2\right)}\in T_{\left(2^2\right);\chi_3,\theta_2},\hbox{\ }
b_{\left(2^2\right)}\notin T_{\left(2^2\right);\chi_3,\theta_2},
$$

$$
a_{\left(2^2\right)}\notin T_{\left(2^2\right);\chi_1,\theta_3},\hbox{\ }
b_{\left(2^2\right)}\notin T_{\left(2^2\right);\chi_1,\theta_3},\hbox{\ }
a_{\left(2^2\right)}\in T_{\left(2^2\right);\chi_2,\theta_3},\hbox{\ }
b_{\left(2^2\right)}\notin T_{\left(2^2\right);\chi_2,\theta_3},
$$
$$
a_{\left(2^2\right)}\notin T_{\left(2^2\right);\chi_3,\theta_3},\hbox{\ }
b_{\left(2^2\right)}\in T_{\left(2^2\right);\chi_3,\theta_3},
$$
so we get the statement about the $(2^2)$-pair.  The statements about the
$(2,1^2)$- and $(1^4)$-pairs are obvious because the stabilizers of any
representatives of these $G$-orbits are trivial.

\enddemo

\proclaim{Corollary 10.1.5} The members of the diameric pairs
$\{a_{\left(2^2\right)},b_{\left(2^2\right)}\}$,
$\{a_{\left(2,1^2\right)},b_{\left(2,1^2\right)}\}$, $\{a_{\left(1^4\right)},
h_{\left(1^4\right)}\}$, $\{b_{\left(1^4\right)},c_{\left(1^4\right)}\}$,
$\{e_{\left(1^4\right)}, f_{\left(1^4\right)}\}$
cannot be distinguished via characters.

\endproclaim

10.2.  The group $G$ of benzene $C_6H_6$ is defined up to conjugation in [7,
IV] (see also [2]) and has order $12$.  We know from the experiment that there
are no chiral pairs among the products of benzene and this conclusion is
confirmed by the Lunn-Senior's model with the fact that there are no groups
$G^\prime$ of order $24$ that contain $G$ (see [7, V]) .  Therefore
$G=G^\prime$ and $Aut_0^\prime(T_{D;G})= Aut_0(T_{D;G})$ for any $D\subset
P_6$.  The classical K\"orner relations among the di- and tri-substitution
homogeneous derivatives (see 7, I) yield immediately

\proclaim{Theorem 10.2.1} If $D=\{(4,2),(3^2)\}$, then the group
$Aut_0^\prime(T_{D;G})$ is the unit group.

\endproclaim

\proclaim{Corollary 10.2.2} The di-substituted (para, ortho, and meta), and
tri-substituted (asymmetric, vicinal, and symmetric) derivatives of benzene are
distinguishable via the K\"orner genetic relations.

\endproclaim

\proclaim{Remark 10.2.3} {\rm In [5], it is shown that the same is true
for any organic compound whose molecule can be divided into a skeleton and six
univalent substituents, and such that it has one mono-substitution and at least
three di-substitution homogeneous derivatives, and, moreover, has a group $G$
of univalent substitution isomerism of order $12$.}

\endproclaim

10.3.  The groups $G$ and $G^\prime$ of univalent substitution isomerism and
stereoisomerism, respectively, of cyclopropane $C_3H_6$ are found up to
conjugation in [7, V], and in [4], and have orders $6$ and $12$, respectively.
Moreover, the group $G$ is dihedral. We set
$$
D=\{(6),(5,1),(4,2),(4,1^2),(3^2)\}.
$$
The possible substitution reactions among the products with empirical formula
in $D$ are represented in the tables below (see [5, Section 5] up to notation).

\vskip 12 pt

\centerline{ Zooms of $T_{D;G}$:}

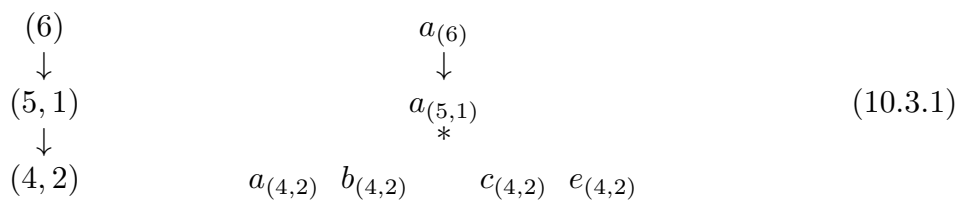
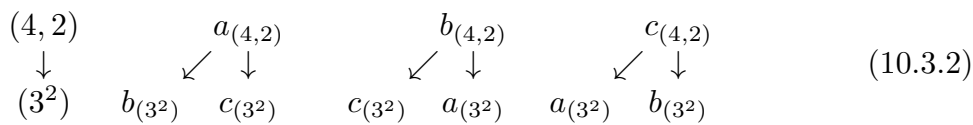
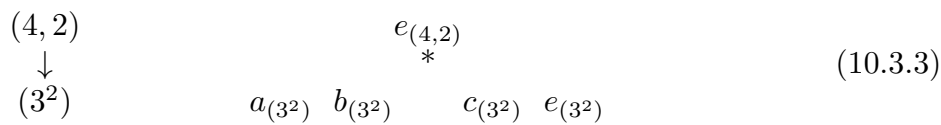
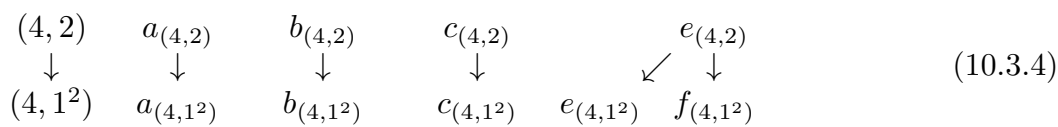
$$
\left.
\matrix \format \c & \c & \c & \c & \c & \c & \c & \c & \c & \c & \c & \c & \c
& \c &  \c & \c &  \c  & \c & \c & \c & \c & \c &  \c  & \c & \c & \c & \c \\
(6) & \hbox{\ \ } & \hbox{\ \ } & \hbox{\ \ } & \hbox{\ \ } & \hbox{\ \ } &
\hbox{\ \ } & \hbox{\ \ } & \hbox{\ \ } & \hbox{\ \ } &
\hbox{\ \ } & \hbox{\ \ } & a_{\left(6\right)} & \hbox{\ \ } & \hbox{\ \ } &
\hbox{\ \ } & \hbox{\ \ } & \hbox{\ \ } & \hbox{\ \ } \cr
\downarrow & \hbox{\ \ } & \hbox{\ \ } & \hbox{\ \ } & \hbox{\ \ } & \hbox{\ \ }
 & \hbox{\ \ } & \hbox{\ \ } & \hbox{\ \ } & \hbox{\ \
} & \hbox{\ \ } & \hbox{\ \ } & \downarrow & \hbox{\ \ } & \hbox{\ \ } &
\hbox{\ \ } & \hbox{\ \ } & \hbox{\ \ } & \hbox{\ \ } \cr
(5,1) & \hbox{\ \ } & \hbox{\ \ } & \hbox{\ \ } & \hbox{\ \ } & \hbox{\ \ } &
\hbox{\ \ } & \hbox{\ \ } & \hbox{\ \ } & \hbox{\ \ } &
\hbox{\ \ } & \hbox{\ \ } & a_{\left(5,1\right)} & \hbox{\ \ } & \hbox{\ \ } &
\hbox{\ \ } & \hbox{\ \ } & \hbox{\ \ } & \hbox{\ \ } \cr
\downarrow & \hbox{\ \ } & \hbox{\ \ } & \hbox{\ \ } & \hbox{\ \ } & \hbox{\ \ }
& \hbox{\ \ } & \hbox{\ \ } & \hbox{\ \ } & \hbox{\ \
} & \hbox{\ \ } & \hbox{ \ \ } & \hbox{ * } & \hbox{\ \ } & \hbox{\ \ } &
\hbox{\ \ } & \hbox{\ \ } & \hbox{\ \ } & \hbox{\ \ } \cr
(4,2) & \hbox{\ \ } & \hbox{\ \ } & \hbox{\ \ } & \hbox{\ \ } & \hbox{\ \ } &
\hbox{\ \ } & \hbox{\ \ } & \hbox{\ \ } & a_{\left(4,2\right)} &\hbox{\ \ } &
b_{\left(4,2\right)} &\hbox{\ } & c_{\left(4,2\right)} & \hbox{\ \ } &
e_{\left(4,2\right)} &\hbox{\ \ } & \hbox{\ \ } & \hbox{\ \ } \cr
\endmatrix \right. \tag 10.3.1
$$

\vskip 15pt

$$
\left.
\matrix \format \c & \c & \c & \c & \c & \c & \c & \c & \c & \c & \c & \c & \c
& \c &  \c & \c &  \c  & \c & \c & \c & \c & \c &  \c  & \c & \c & \c & \c \\
(4,2) & \hbox{\ \ } & \hbox{\ \ } & \hbox{\ \ } & \hbox{\ \ } &
a_{\left(4,2\right)} &
\hbox{\ \ } & \hbox{\ \ } & \hbox{\ \ } & \hbox{\ \ } & \hbox{\ \ } &
b_{\left(4,2\right)}
& \hbox{\ \ } & \hbox{\ \ } & \hbox{\ \ } & \hbox{\ \ } &
c_{\left(4,2\right)} & \hbox{\ \ }  \cr
\downarrow & \hbox{\ \ } & \hbox{\ \ } & \hbox{\ \ } & \swarrow & \downarrow &
\hbox{\ \ } & \hbox{\ \ } & \hbox{\ \ } & \hbox{\ \ } & \swarrow & \downarrow &
\hbox{\ \ } & \hbox{\ \ } & \hbox{\ \ } & \swarrow
& \downarrow & \hbox{\ \ } \cr
(3^2) & \hbox{\ \ } & \hbox{\ \ } & b_{\left(3^2\right)} & \hbox{\ \ } &
c_{\left(3^2\right)} & \hbox{\ \ } & \hbox{\ \ } & \hbox{\ \ } &
c_{\left(3^2\right)}
& \hbox{\ \ } & a_{\left(3^2\right)} & \hbox{\ \ }
& \hbox{\ \ } & a_{\left(3^2\right)} & \hbox{\ \ } &
b_{\left(3^2\right)} & \hbox{\ \ }    \cr
\endmatrix \right.\tag 10.3.2
$$

\vskip 12pt

$$
\left.
\matrix \format \c & \c & \c & \c & \c & \c & \c & \c & \c & \c & \c & \c & \c
& \c &  \c & \c &  \c  & \c & \c & \c & \c & \c &  \c  & \c & \c & \c & \c \\
(4,2) & \hbox{\ \ } & \hbox{\ \ } & \hbox{\ \ } & \hbox{\ \ } & \hbox{\ \ } &
\hbox{\ \ } & \hbox{\ \ } & \hbox{\ \ } & \hbox{\ \ } &
\hbox{\ \ } & \hbox{\ \ } & e_{\left(4,2\right)} & \hbox{\ \ } & \hbox{\ \ } &
\hbox{\ \ } & \hbox{\ \ } & \hbox{\ \ } & \hbox{\ \ } \cr
\downarrow & \hbox{\ \ } & \hbox{\ \ } & \hbox{\ \ } & \hbox{\ \ } & \hbox{\ \ }
& \hbox{\ \ } & \hbox{\ \ } & \hbox{\ \ } & \hbox{\ \
} & \hbox{\ \ } & \hbox{ \ \ } & \hbox{ * } & \hbox{\ \ } & \hbox{\ \ } &
\hbox{\ \ } & \hbox{\ \ } & \hbox{\ \ } & \hbox{\ \ } \cr
(3^2) & \hbox{\ \ } & \hbox{\ \ } & \hbox{\ \ } & \hbox{\ \ } & \hbox{\ \ } &
\hbox{\ \ } & \hbox{\ \ } & \hbox{\ \ } & a_{\left(3^2\right)} &\hbox{\ \ } &
b_{\left(3^2\right)} &\hbox{\ } & c_{\left(3^2\right)} & \hbox{\ \ } &
e_{\left(3^2\right)} &\hbox{\ \ } & \hbox{\ \ } & \hbox{\ \ } \cr
\endmatrix \right.\tag 10.3.3
$$

\vskip 12pt

$$
\left.
\matrix \format \c & \c & \c & \c & \c & \c & \c & \c & \c & \c & \c & \c & \c
& \c &  \c & \c &  \c  & \c & \c & \c & \c & \c &  \c  & \c & \c & \c & \c \\
(4,2) & \hbox{\ \ } &
\hbox{\ \ } & a_{\left(4,2\right)} &
\hbox{\ \ } & \hbox{\ \ } & \hbox{\ \ } & b_{\left(4,2\right)} & \hbox{\ \ } &
\hbox{\ \ } &\hbox{\ \ } & c_{\left(4,2\right)} & \hbox{\ \ } &
\hbox{\ \ } & \hbox{\ \ } & \hbox{\ \ } & e_{\left(4,2\right)} &
\hbox{\ \ } & \hbox{\ \ }
\cr
\downarrow & \hbox{\ \ } &
\hbox{\ \ } & \downarrow & \hbox{\ \ } & \hbox{\ \
} &\hbox{\ \ } & \downarrow & \hbox{\ \ } & \hbox{\ \ } &\hbox{\ \ }
&\downarrow &\hbox{\ \ } &
\hbox{\ \ } & \hbox{\ \ } & \swarrow & \downarrow & \hbox{\ \ } &
\hbox{\ \ }
\cr
(4,1^2) & \hbox{\ \ } &
\hbox{\ \ } & a_{\left(4,1^2\right)} & \hbox{\ \ } & \hbox{\ \ } &
\hbox{\ \ } & b_{\left(4,1^2\right)} & \hbox{\ \ } & \hbox{\ \ } &
\hbox{ \ \ } & c_{\left(4,1^2\right)} &\hbox{\ \ } &
\hbox{\ \ } & e_{\left(4,1^2\right)} & \hbox{\ \ } &
f_{\left(4,1^2\right)} & \hbox{\ \ } & \hbox{\ \ }
\cr
\endmatrix \right.\tag 10.3.4
$$

\vskip 12pt

The $G^\prime$-orbits are
$$
a_{\left(4,2\right)}\cup b_{\left(4,2\right)},\hbox{\ }
c_{\left(4,2\right)}, \hbox{\ } e_{\left(4,2\right)},
$$
$$ a_{\left(3^2\right)}\cup
b_{\left(3^2\right)}, \hbox{\ } c_{\left(3^2\right)}, \hbox{\ }
e_{\left(3^2\right)},
$$
$$ a_{\left(4,1^2\right)}\cup b_{\left(4,1^2\right)},
\hbox{\ } c_{\left(4,1^2\right)}, \hbox{\ } e_{\left(4,1^2\right)}\cup
f_{\left(4,1^2\right)}.
$$
In particular, the pairs
$\{a_{\left(4,2\right)}, b_{\left(4,2\right)}\}$,
$\{a_{\left(3^2\right)}, b_{\left(3^2\right)}\}$,
$\{a_{\left(4,1^2\right)}, b_{\left(4,1^2\right)}\}$,
$\{e_{\left(4,1^2\right)}, f_{\left(4,1^2\right)}\}$,
are chiral.

\proclaim{Theorem 10.3.5} (i) One has
$$
Aut_0(T_{D;G})\simeq S_3\times S_2,
$$
where the permutations of $S_3$ simultaneously move the elements $a$, $b$, $c$
with indices $(4,2)$, $(3^2)$, and $(4,1^2)$, and $S_2$ consists of all
permutations of $e_{\left(4,1^2\right)}, f_{\left(4,1^2\right)}$;

(ii) one has
$$
Aut_0^\prime(T_{D;G})\simeq S_2\times S_2,
$$
where the permutations of the first direct component $S_2$ simultaneously move
the elements $a$, $b$ with indices $(4,2)$, $(3^2)$, and $(4,1^2)$, and the
second direct component $S_2$ consists of all permutations of
$e_{\left(4,1^2\right)}, f_{\left(4,1^2\right)}$.

\endproclaim

\demo{Proof} Let $\alpha\in Aut_0(T_{D;G})$. In accordance with (10.3.1),
$\alpha(a_{\left(6\right)})=a_{\left(6\right)}$,
$\alpha(a_{\left(5,1\right)})=a_{\left(5,1\right)}$, and the diagrams (10.3.2)
and (10.3.3) yield
$\alpha(e_{\left(3^2\right)})=e_{\left(3^2\right)}$ and
$\alpha(e_{\left(4,2\right)})=e_{\left(4,2\right)}$.
The
last equality together with (10.3.4) implies that the sets
$\{e_{\left(4,1^2\right)}, f_{\left(4,1^2\right)}\}$, and
$\{a_{\left(4,1^2\right)}, b_{\left(4,1^2\right)}, c_{\left(4,1^2\right)}\}$
are $\alpha$-stable. Now, comparing (10.3.2) with (10.3.4), we get (i).
Part (ii) can be obtained from (i) by taking into account that a chiral
automorphism $\alpha\in Aut_0(T_{D;G})$ maps a chiral pair onto a chiral pair.

\enddemo

\proclaim{Corollary 10.3.6} The products that correspond to the different sets
of structural formulae below are distinguishable via substitution reactions
among the elements of $T_{D;G}$:
$$
\{a_{\left(4,2\right)}, b_{\left(4,2\right)}\},\hbox{\ }
\{c_{\left(4,2\right)}\},\hbox{\ } \{e_{\left(4,2\right)}\},
$$
$$
\{a_{\left(3^2\right)},b_{\left(3^2\right)}\}, \hbox{\ }
\{c_{\left(3^2\right)}\},\hbox{\ }\{e_{\left(3^2\right)}\},
$$
$$
\{a_{\left(4,1^2\right)},
b_{\left(4,1^2\right)}\}\hbox{\ }
\{c_{\left(4,1^2\right)}\}, \hbox{\ }
\{e_{\left(4,1^2\right)},
f_{\left(4,1^2\right)}\}.
$$
The products that correspond to the members of a particular set are
indistinguishable via substitution reactions.

\endproclaim

\demo{Proof} For the first statement, it is enough to note that all sets are
$Aut_0^\prime(T_{D;G})$-stable.  Since all pairs are chiral, we get the second
statement by using Theorem 8.2.6.

\enddemo

Theorem 2.1 and Section 5 from [5], and [4, 7.3.1] yield immediately

\proclaim{Remark 10.3.7} {\rm All statements in this Subsection are true for
any organic compound whose molecule can be divided into a skeleton and six
univalent substituents, and such that it has one mono-substitution and at least
three di-substitution homogeneous derivatives, and, moreover, has a group
$G\leq S_6$ of univalent substitution isomerism of order $6$, that is dihedral.
}

\endproclaim

\proclaim{Remark 10.3.8} {\rm The more precise definition of distinguishability
of two compounds that includes chirality, Corollary 10.3.6, and Remark 10.3.7,
specify the conclusions in [5, Section 5].}

\endproclaim

In [7, VII], the group $G^{\prime\prime}$ of structural isomerism of
cyclopropane is described.  In particular, the order of $G^{\prime\prime}$ is
$48$, and we have $G\leq G^\prime\leq G^{\prime\prime}$.  The
$G^{\prime\prime}$-orbits in $T_D$ are
$$
a_{\left(6\right)},
$$
$$
a_{\left(5,1\right)},
$$
$$
a_{\left(4,2\right)}\cup
b_{\left(4,2\right)}\cup
e_{\left(4,2\right)},\hbox{\ }
c_{\left(4,2\right)},
$$
$$
a_{\left(3^2\right)}\cup
b_{\left(3^2\right)}, \hbox{\ }
c_{\left(3^2\right)}\cup
e_{\left(3^2\right)},
$$
$$
a_{\left(4,1^2\right)}\cup
b_{\left(4,1^2\right)}\cup
e_{\left(4,1^2\right)}\cup
f_{\left(4,1^2\right)},\hbox{\ }
c_{\left(4,1^2\right)}.
$$
Therefore we obtain the next proposition which completes the description of the
isomers of cyclopropane, considered in present Subsection.

\proclaim{Proposition 10.3.9} The products of cyclopropane that correspond to
the structural formulae within any one of the sets
$$
\{a_{\left(6\right)}\},
$$
$$
\{a_{\left(5,1\right)}\},
$$
$$
\{a_{\left(4,2\right)},
b_{\left(4,2\right)},
e_{\left(4,2\right)}\},\hbox{\ }
\{c_{\left(4,2\right)}\},
$$
$$
\{a_{\left(3^2\right)},
b_{\left(3^2\right)}\}, \hbox{\ }
\{c_{\left(3^2\right)},
e_{\left(3^2\right)}\},
$$
$$
\{a_{\left(4,1^2\right)},
b_{\left(4,1^2\right)},
e_{\left(4,1^2\right)},
f_{\left(4,1^2\right)}\},\hbox{\ }
\{c_{\left(4,1^2\right)}\}
$$
are structurally identical, and the products that correspond to members of
different sets and have the same empirical formula, are structural isomers.

\endproclaim

\heading {References}

\endheading

\noindent [1] P. Doubilet, G.-C. Rota, R. Stanley, On the foundations of
combinatorial theory (VI): the idea of generating function, {\it in}
Proceedings of the Sixth Berkeley Symposium on Mathematical Statistics and
Probability, vol. II, University of California Press, Berkeley and Los Angeles,
1972.

\noindent [2] W. H\"asselbarth, The Inverse Problem of Isomer Enumeration, J.
Comput. Chem. 8 (1987), 700 -- 717.

\noindent [3] V.  V.  Iliev, On Lunn-Senior's Mathematical Model of Isomerism
in Organic Chemistry.  Part I, MATCH - Commun.  Math.  Comput.  Chem.  40
(1999), 153 -- 186.

\noindent [4] V.  V.  Iliev, On the Inverse Problem of Isomer Enumeration:
Part II, Case of Cyclopropane, MATCH - Commun.  Math.  Comput.  Chem.  43
(2001), 79 -- 84.

\noindent [5] V.  V.  Iliev, On certain organic compounds with one
mono-substitution and at least three di-substitution homogeneous derivatives,
J. Math. Chem. (to appear).

\noindent [6]  A.  Kerber, Applied Finite Group Actions, Springer, Berlin,
1999.

\noindent [7] A. C. Lunn, J. K. Senior, Isomerism and Configuration, J.
Phys. Chem. 33 (1929), 1027 - 1079.

\noindent [8] I.  G.  Macdonald, Symmetric Functions and Hall Polynomials,
Clarendon Press, Oxford, 1995.

\end